\documentclass[11pt]{amsart}
\usepackage{amssymb,latexsym}%\documentclass[12pt,leqno]{amsart}

\def\triangle{\Delta}

\def\be1{{\begin{equation}}}
\def\ee1{{\end{equation}}}

\def\part{\partial}

\def\ba{\begin{array}}
\def\ea{\end{array}}

\newtheorem{corollary}{Corollary}[section]

\numberwithin{equation}{section}

\newtheorem{lemma}{Lemma}[section]
\newtheorem{proposition}[lemma]{Proposition}
\newtheorem{theorem}[lemma]{Theorem}

\title[GW conformal invariant for closed four dimensional submanifolds]
{Graham-Witten's conformal invariant for closed four dimensional submanifolds}

\author{Yongbing Zhang}
\address{School of Mathematical Sciences and Wu Wen-Tsun Key Laboratory of Mathematics,
University of Science and Technology of China, Hefei 230026, Anhui Province, China}
\email{ybzhang@amss.ac.cn}
\thanks{The project is supported by ¡°the Fundamental Research Funds for
the Central Universities¡±}


\begin{document}
\maketitle

{\sl Dedicated to Professor A. Chang and P. Yang on the occasion of their 70th birthdays.}

\begin{abstract}
It was proved by Graham and Witten in 1999 that conformal invariants of submanifolds can be obtained
via volume renormalization of  minimal surfaces in conformally compact Einstein manifolds.
The conformal invariant of a submanifold $\Sigma$ is contained in the volume expansion of the minimal surface
which is asymptotic to $\Sigma$ when the minimal surface approaches the conformaly infinity.
In the paper we give the explicit expression of Graham-Witten's conformal invariant for closed four dimensional submanifolds and
find critical points of the conformal invariant in the case of Euclidean ambient spaces.
\end{abstract}

\section{Introduction}

In the introduction we give a description of the main result and some related background of the paper.
The terminologies used in the introduction will be recalled in the next section.

Let $(X^{d+1},g_+)$ be a conformally compact Einstein manifold and $(M^d,[g_{confinf}])$ its conformal infinity.
A given metric $\overline{g}\in [g_{confinf}]$ uniquely determines a special defining function $r$ on a neighborhood of $M$ in $\overline{X}$,
upon to the conditions that $(r^2g_+)|_{M}=\overline{g}$ and $|dr|_{r^2g_+}=1$ \cite{Graham00}.
We denote
%\begin{equation}\label{gc}
$g_c=r^2g_+$.
%\end{equation}
With the special defining function $r$, one can identify $M\times [0,\epsilon)$, for some $\epsilon>0$,
with a neighborhood of $M$  in $\overline{X}$. We denote the neighborhood by $X_\epsilon$, and the identification
\begin{equation}\label{identification}
M\times [0,\epsilon)\cong X_\epsilon
\end{equation}
is defined as follows: $(p, r)\in M\times [0,\epsilon)$ corresponds to the point obtained by following the flow of $\nabla^{g_c}r$
emanating from $p$ for $r$ units of time.
$g_c$ on $M\times [0,\epsilon)$ takes the form of
\begin{equation}\label{gr}
g_c=dr^2+\widetilde{g},
\end{equation}
where $\widetilde{g}$ is a $1$-parameter family of metrics on $M$ with the parameter $r$.
By solving the Einstein equation $Ric(g_+)=-dg_+$, for $d$ odd the expansion of $\widetilde{g}$ is of the form
\begin{equation}\label{oddexpansion}
\widetilde{g}=g^{(0)}+g^{(2)}r^2+(\mbox{even powers})+g^{(d-1)}r^{d-1}+g^{(d)}r^d+\cdots
\end{equation}
where $g^{(j)}$ are tensors on $M$ and the dots stand for terms vanishing to higher order.
For $j$ even and $0\leq j\leq d-1$, the tensor $g^{(j)}$
is locally formally determined by the boundary value $g^{(0)}=\overline{g}$,
but $g^{(d)}$ is formally undetermined; for $d$ even the expansion is
\begin{equation}\label{evenexpansion}
\widetilde{g}=g^{(0)}+g^{(2)}r^2+(\mbox{even powers})+hr^d\log r+g^{(d)}r^{d}+\cdots
\end{equation}
where $g^{(j)}$ and $h$ are locally formally determined for $j$ even and $0\leq j\leq d-2$ by $g^{(0)}=\overline{g}$.

The main object of the paper is a minimal surface in the conformally compact Einstein manifold $(X,g_+)$
with prescribed asymptotic boundary.
Let $\Sigma^n$ be a submanifold of $M$ and $Y^{n+1}\hookrightarrow (X,g_+)$ be a minimal surface which is asymptotic to $\Sigma$.
The problem of existence and regularity of such minimal surfaces has been studied by Anderson \cite{Anderson82, Anderson83},
Hardt-Lin \cite{HardtLin}, Lin \cite{Lin89CPAM, Lin89Invent, Lin2012}, Tonegawa \cite{Tonegawa}, Han-Jiang \cite{HanJiang} and
Han-Shen-Wang \cite{HanShenWang}.
We denote
\begin{equation}\label{g}
g=\overline{g}|_\Sigma.
\end{equation}
The connections with respect to $(M,\overline{g})$ and $(\Sigma, g)$ will be denoted by $\overline{\nabla}$ and $\nabla$ respectively,
and the connection of the normal bundle $T^\perp \Sigma$ of the immersion $\Sigma\hookrightarrow (M^d, \overline{g})$ will be denoted by $\nabla^\perp$.

Graham and Witten \cite{GrahamWitten} have introduced a natural and useful way to reformulate $Y$. Namely,
near the boundary $M$ they express $Y$ as a graph over $\Sigma\times [0,\epsilon)$ and expand the height functions of the graph
in $r$. Near a point of $\Sigma^n$, let $(x^i, y^\alpha)$ be a local coordinate chart of $M^d$, where $1\leq i\leq n$ and $n+1\leq \alpha\leq d$, so that
\begin{equation}\label{coordinatexy}
\Sigma=\{y=0\}; \quad \overline{g}(\frac{\partial}{\partial x^i},\frac{\partial}{\partial y^\alpha})=0\,\ \mbox{ on}\,\ \Sigma\,\ ,\forall \,\ i, \alpha.
\end{equation}
Note that via the identification (\ref{identification}),
one has an extension of the coordinates $(x^i,y^\alpha)$ into $X$, which together with $r$ forms
a local coordinate chart of $\overline{X}$.
The minimal surface $Y$ can be written as a graph $\{y^\alpha=u^\alpha(x,r)\}$.
That is, near the boundary $Y=(x^i,u^\alpha(x,r),r)$.

Graham and Witten \cite{GrahamWitten} proved that for $n$ odd
\begin{equation}\label{uexpansionoddn}
u=u^{(2)}r^2+(\mbox{even powers})+u^{(n+1)}r^{n+1}+u^{(n+2)}r^{n+2}+\cdots,
\end{equation}
and for $n$ even
\begin{equation}\label{uexpansionevenn}
u=u^{(2)}r^2+(\mbox{even powers})+u^{(n)}r^{n}+w_nr^{n+2}\log r+u^{(n+2)}r^{n+2}+\cdots,
\end{equation}
where the $u^{(k)}=(u^{(k)\alpha}), k<n+2$, and $w_n=(w_n^\alpha)$ are functions of $\Sigma$ and locally determined,
while $u^{(n+2)}$ is not locally determined.
They showed that the first non-vanishing coefficient $u^{(2)}$ is related to the mean curvature of $\Sigma\hookrightarrow (M,\overline{g})$ by
\begin{equation}
u^{(2)}=\frac{1}{2n}H.
\end{equation}

In the paper, we calculate the coefficient $u^{(4)}$ for $n\geq 3$ by using the graphic minimal surface equation of $Y$.
\begin{theorem}\label{u4prop}
For $n\geq 3$,
\begin{eqnarray}\label{u4intro}
8n(n-2)u^{(4)}&=&\triangle^\perp H^\alpha+g^{ij}g^{kl}\langle A_{ik},H\rangle A_{jl}^\alpha-\frac{2}{n^2}|H|^2H^\alpha \nonumber
\\&&+g^{ij}\overline{g}^{\alpha\gamma}H^\beta\overline{W}_{i\beta j\gamma}+(n-4)\overline{g}^{\alpha\gamma}\overline{P}_{\beta\gamma}H^\beta
-g^{ij}\overline{P}_{ij}H^\alpha+2ng^{ik}g^{jl}\overline{P}_{kl}A_{ij}^\alpha
\\&&+n\overline{g}^{\alpha\gamma}g^{ij}(\overline{\nabla}_\gamma \overline{P}_{ij}-2\overline{\nabla}_i\overline{P}_{j\gamma})
-\frac{n-2}{n}H^\beta H^\gamma \overline{\Gamma}^\alpha_{\beta\gamma},\nonumber
\end{eqnarray}
where $A_{ij}$ denotes the second fundamental form of $\Sigma\hookrightarrow (M,\overline{g})$,
$\overline{\Gamma}, \overline{P}$ and $\overline{W}$ denote the Christoffel symbol, Schouten and Weyl tensor of $\overline{g}$ respectively.
\end{theorem}

For $n=2$, the calculation of $w_2$ in (\ref{uexpansionevenn}) can be carried out in a similar way. We have
\begin{equation}\label{v}
w_2=-\frac{1}{16}W,
\end{equation}
where
\begin{eqnarray}\label{W}
W^\alpha&=&\triangle^\perp H^\alpha+g^{ij}g^{kl}\langle A_{ik},H\rangle A_{jl}^\alpha-\frac{1}{2}|H|^2H^\alpha \nonumber
\\&&+g^{ij}\overline{g}^{\alpha\gamma}H^\beta\overline{W}_{i\beta j\gamma}-2\overline{g}^{\alpha\gamma}\overline{P}_{\beta\gamma}H^\beta
-g^{ij}\overline{P}_{ij}H^\alpha+4g^{ik}g^{jl}\overline{P}_{kl}A_{ij}^\alpha
\\&&+2\overline{g}^{\alpha\gamma}g^{ij}(\overline{\nabla}_\gamma \overline{P}_{ij}-2\overline{\nabla}_i\overline{P}_{j\gamma}). \nonumber
\end{eqnarray}
Note that if $(X, g_+)$ is the Poincar\'{e} half-plane model of the hyperbolic space $\mathbb{H}^{d+1}$, for which $(M,\overline{g})=\mathbb{R}^d$
and $g^{(2)}=-\overline{P}=0$, then for $\Sigma^2$ we have
$$W=\triangle^\perp H+g^{ij}g^{kl}\langle A_{ik},H\rangle A_{jl}-\frac{1}{2}|H|^2H.$$
A surface $\Sigma^2\hookrightarrow \mathbb{R}^d$ with $W=0$ is well-known to be a Willmore surface.
Han and Jiang \cite{HanJiang} proved that in the case of $\Sigma^2\hookrightarrow \mathbb{R}^3$
the log term vanishes if and only if $\Sigma$ is a Willmore surface.

Our main aim of the paper is to calculate the conformal invariant,
introduced by Graham and Witten (for each dimensional submanifold) \cite{GrahamWitten},
for closed four dimensional submanifolds. Here the conformal invariant for submanifolds is a functional of submanifolds which is invariant under
conformal transformations of the metric of the ambient space. Recently, conformal invariants for hypersurfaces, constructed
by using the volume renormalization of solutions to the singular Yamabe problem or general singular volume measures, and related aspects
have been extensively studied \cite{GlarosGoverHalbachWaldron, Graham17, GoverWaldron1407, GoverWaldron1506, GoverWaldron1603,
GoverWaldron1611, GoverWaldron1611B, Vyatkin}.
Graham-Witten's conformal invariants are obtained from the renormalization process of the volume
of a minimal surface $Y^{n+1}$ in a conformally compact Einstein manifold $(X^{d+1}, g_+)$ which is asymptotic to a submanifold
$\Sigma^n$ immersed in the conformal infinity of $X$.

Near the asymptotic boundary, the volume form of $Y$ takes the form of
\begin{equation}\label{volumeelementexpansion}
d\mu_Y=r^{-n-1}[v^{(0)}+v^{(2)}r^2+(\mbox{even powers})+v^{(n)}r^n+\cdots]d\mu_{\Sigma} dr,
\end{equation}
where $v^{(j)}$ are locally determined functions of $\Sigma$, $v^{(n)}=0$ for $n$ odd, and $d\mu_\Sigma$ is the volume form of $(\Sigma, g)$.
As $\epsilon\rightarrow 0$ and for $n$ odd, the volume
\begin{equation}\label{areaexpansionodd}
\mbox{Vol}_{g_+}(Y\cap \{r>\epsilon\})=c_0\epsilon^{-n}+c_2\epsilon^{-n+2}+(\mbox{odd powers})+c_{n-1}\epsilon^{-1}+c_n+o(1),
\end{equation}
and for $n$ even
\begin{equation}\label{areaexansioineven}
\mbox{Vol}_{g_+}(Y\cap \{r>\epsilon\})=c_0\epsilon^{-n}+(\mbox{even powers})+c_{n-2}\epsilon^{-2}+L_n\log\frac{1}{\epsilon}+c_n+o(1),
\end{equation}
where
\begin{equation}
L_n=\int_\Sigma v^{(n)}d\mu_\Sigma.
\end{equation}
An area law was proposed by Ryu and Takayanagi \cite{RyuTakayanagiPRL06, RyuTakayanagiJHEP06} which holographically identifies
the volume in (\ref{areaexpansionodd}) or (\ref{areaexansioineven}) with the entanglement entropy in quantum (conformal) field theories.
Graham and Witten proved the following
\begin{theorem} (\cite{GrahamWitten}) \label{conformalinvariants}
If $n$ is odd, then $c_n$ is independent of the choice of special defining function.
If $n$ is even, then $L_n$ is independent of the choice of special defining function.
\end{theorem}

Note that there is the one-to-one correspondence between representative metrics of the conformal class $(M, [g_{confinf}])$ and special defining functions.
Hence $c_n$ for $n$ odd and $L_n$ for $n$ even are conformal invariants. $c_n$ is called the renormalized volume of $Y$.
For $n=1$, an explicit expression of the renormalized volume $c_1$ has been obtained by Alexakis and Mazzeo \cite{AlexakisMazzeo}.
For $n=2$, Graham and Witten \cite{GrahamWitten} showed that
$$-8L_2=\int_\Sigma(|H|^2+4g^{ij}\overline{P}_{ij})d\mu_\Sigma,$$
where $\overline{P}=\frac{1}{d-2}(\overline{Ric}-\frac{1}{2(d-1)}\overline{R}\overline{g})$ is the Schouten tensor of $(M^d, \overline{g})$.
$L_2$ is closely related to the Willmore functional.
For example, when $\Sigma$ is a closed two dimensional hypersurface in $(M^3, \overline{g})$,
it easily follows from the Gauss formula and Gauss-Bonnet formula that
$$-8L_2=8\pi \chi(\Sigma)+2\int_\Sigma|\stackrel{\circ}{A}|^2d\mu_\Sigma,$$
where $\stackrel{\circ}{A}$ denotes the traceless part of the second fundamental form of $\Sigma$.

Volume renormalization of a conformally compact Einstein manifold had been studied at almost the same time
\cite{HenningsonSkenderis9806, HenningsonSkenderis9812, Graham00},
which can be viewed as the extreme case under the setting of the paper: $\Sigma=M$ and $Y=X$.
The log terms of the volume expansion were calculated explicitly in lower dimensions. For example, in dimension two
$$L_{d=2}=-\frac{1}{4}\int_{M^2}R_{\overline{g}}d\mu_{\overline{g}}%=-\pi \chi(M^2)
,$$
and in dimension four
\begin{equation}L_{d=4}=\frac{1}{4}\int_{M^4}\sigma_2(\overline{P})d\mu_{\overline{g}} \label{Ld4}
%=\frac{1}{2}\pi^2\chi(M^4)-\frac{1}{64}\int_{M^4}|\overline{W}|^2d\mu_{\overline{g}}
,\end{equation}
where $\overline{P}$ is the Schouten tensor of $\overline{g}$ and $\sigma_2$ is the second
elementary symmetric function.
The geometry and topology of closed four dimensional manifolds which admit a metric $\overline{g}$ such that
$\int_{M^4} \overline{R}d\mu_{\overline{g}}>0$ and $\int_{M^4} \sigma_2(\overline{P})d\mu_{\overline{g}}>0$
has been studied by Gursky \cite{Gursky} and Chang-Gursky-Yang \cite{ChangGurskyYang}.
For discussions on the renormalized volume of conformally compact Einstein manifold, see for instance
\cite{Albin, Anderson01, ChangFangGraham, YangKingChang, FeffermanGraham02, GrahamZworski}.
We believe that the renormalized volume corresponds to the log-determinant of a pseudo-differential operator of orders $d$
at the quantum field theory side. A likely choice of the pseudo-differential operator could be the linearization of
the Dirichlet-to-Neumann operator given by the conformally compact Einstein metric.

In the paper, we calculate $L_4$ for closed four dimensional submanifolds.

\begin{theorem}\label{GWinvariant4}
Let $\Sigma^4$ be a closed submanifold in $(M^d, \overline{g})$, we have
\begin{eqnarray}
128L_4&=&\int_\Sigma (|\nabla^\perp H|^2-g^{ik}g^{jl}\langle A_{ij},H\rangle \langle A_{kl},H\rangle +\frac{7}{16}|H|^4)d\mu_g \nonumber
\\&&+16\int_\Sigma [(g^{ij}\overline{P}_{ij})^2-g^{ik}g^{jl}\overline{P}_{ij}\overline{P}_{kl}
+g^{ij}\overline{g}^{\alpha\beta}\overline{P}_{i\alpha}\overline{P}_{j\beta}-\frac{1}{d-4}g^{ij}\overline{B}_{ij}]d\mu_g \label{L4full}
\\&&+\int_\Sigma (-16 g^{ik}g^{jl}\overline{P}_{ij}\langle A_{kl},H \rangle+5g^{ij}\overline{P}_{ij}|H|^2+8\overline{P}(H,H)-g^{ij}\overline{W}_{i\alpha j\beta}H^\alpha H^\beta)d\mu_g \nonumber
\\&&-8\int_\Sigma g^{ij}H^\beta(\overline{\nabla}_\beta \overline{P}_{ij}-2\overline{\nabla}_j\overline{P}_{i\beta})d\mu_g,\nonumber
\end{eqnarray}
where $\overline{B}_{ij}$ is the Bach tensor of $(M^d, \overline{g})$.
\end{theorem}

When $(X^{d+1}, g_+)$ is the Poincar\'{e}  half-plane model with $(M^d, \overline{g})=\mathbb{R}^d$, $L_4$ is simplified to
\begin{equation}\label{pureextrinsic}
\frac{1}{128}\int_\Sigma (|\nabla^\perp H|^2-g^{ik}g^{jl}\langle A_{ij},H\rangle \langle A_{kl},H\rangle +\frac{7}{16}|H|^4)d\mu_g.
\end{equation}
In the last part of the paper, we will consider the functional
\begin{equation}\label{L4functional}
\mathcal{L}_4(\Sigma^4):=\frac{1}{2}\int_\Sigma (|\nabla^\perp H|^2-g^{ik}g^{jl}\langle A_{ij},H\rangle \langle A_{kl},H\rangle +\frac{7}{16}|H|^4)d\mu_g,
\end{equation}
for closed four dimensional submanifolds $\Sigma^4\hookrightarrow \mathbb{R}^d$.
Guven \cite{Guven} has constructed a bending energy for closed four dimensional submanifolds immersed in
$\mathbb{R}^5$, which is invariant under special conformal transformations of $\mathbb{R}^5$ and reads
\begin{equation}\label{Guvenfunctional}
H_2=\frac{1}{2}\int_\Sigma (|\nabla  H|^2-|A|^2H^2+\frac{7}{16}|H|^4)d\mu_g.
\end{equation}
(A factor of $-2$ had been dropped in the original paper \cite{Guven}, pointed out by Graham and Reichert \cite{GrahamReichert}.)
Note also that Graham and Reichert \cite{GrahamReichert} prove that up to a multiple
the gradient of $L_n$ is given by $w_n$ in (\ref{uexpansionevenn}).

In general, the functional $\mathcal{L}_4$ is unbounded from below and above even for closed four dimensional submanifolds
with a fixed topology. Let $\mathbb{S}^k(r)$ denote the round sphere in $\mathbb{R}^{k+1}$ of radius $r$.
Examples of critical points of $\mathcal{L}_4$ include
$\mathbb{S}^4$,
$\mathbb{S}^3(1)\times \mathbb{S}^1(\frac{1}{\sqrt{3}})$, $\mathbb{S}^3(1)\times \mathbb{S}^1(\sqrt{\frac{3}{5}})$,
$\mathbb{S}^2(1)\times \mathbb{S}^2(1)$,
$\mathbb{S}^2(1)\times \mathbb{S}^1(\frac{1}{\sqrt{2}})\times \mathbb{S}^1(\frac{1}{\sqrt{2}})$,
$\mathbb{S}^2(1)\times \mathbb{S}^1(\frac{1}{\sqrt{2}})\times \mathbb{S}^1(\frac{3}{\sqrt{10}})$,
$\mathbb{S}^1(1)\times \mathbb{S}^1(1)\times \mathbb{S}^1(1)\times \mathbb{S}^1(1)$, and
$\mathbb{S}^1(1)\times \mathbb{S}^1(1)\times \mathbb{S}^1(1)\times \mathbb{S}^1(\frac{3}{\sqrt{5}})$.
These are all the closed and four dimensional critical points of $\mathcal{L}_4$ which is a product of round spheres.

The calculations in the paper are elementary.
In Section 2, we give a brief recall of conformally compact Einstein manifolds and minimal surfaces in
conformally compact Einstein manifolds. We identify the constituent components of $u^{(4)}$ for $n\geq 3$ by using the graphic
minimal surface equation. In Section 3, we reformulate $u^{(4)}$ in a covariant form and prove Proposition \ref{u4prop}.
In Section 4, we calculate the coefficient $v^{(4)}$ in (\ref{volumeelementexpansion}) for $n\geq 4$, from which Proposition
\ref{GWinvariant4} follows so that we obtain Graham-Witten's conformal invariant  $L_4$ for
closed four dimensional submanfolds. In the last section, we consider Graham-Witten's conformal invariant for closed four dimensional
submanifolds of Euclidean spaces and find simple critical points of it.

{\bf Acknowledgements.} The author would like to thank Professor Qing Han
for his stimulating talk which also brought the work \cite{GrahamWitten} of Graham and Witten to his attention.
This work is completed when the author visits Princeton University. He would like
to thank the Department of Mathematics for its hospitality.

{\bf Notice.} After the completion of the manuscript, we learn that, along with other interesting results and discussions,
$u^{(4)}$ and $L_4$ have also been calculated by Graham and Reichert \cite{GrahamReichert}. The critical points
contained in the paper have also been found in \cite{GrahamReichert}.

\section{Minimal submanifolds in conformally compact Einstein manifolds}

In the first part of this section, referring mainly to \cite{Graham00, GrahamWitten},
we give a brief review to the related aspects of conformally compact Einstein manifolds,
and minimal submanifolds in conformally compact Einstein manifolds.
In the second part, we identify the constituent components of $u^{(4)}$ from the graphic minimal surface equation of $Y$.

\subsection{Minimal surfaces in a conformally compact Einstein manifold}
Let $\overline{X}$ be an $(d+1)$-dimensional manifold with a $d$-dimensional boundary $M$.
We call a function $r\in C^\infty(\overline{X})$ a defining function for
$M$ if it satisfies that $r|_X>0$, $r|_M=0$ and $(dr)|_M\neq 0$.
A metric $g_+$ defined on $X$ is said to be conformally compact if
$r^2g_+$ extends to a metric on $\overline{X}$.
Given a conformally compact manifold $(X,g_+)$, $(r^2g_+)|_{TM}$ induces a conformal class of metrics as
$r$ runs over the space of defining functions, called the conformal infinity of $(X,g_+)$ and denoted by $(M, [g_{confinf}])$.
Let $g_c=r^2g_+$.
The sectional curvature of $(X, g_+)$ is asymptotic to $-(|dr|^2_{g_c})|_M$,
which is independent of the the choices of $r$, as $r\rightarrow 0$.

If a conformally compact manifold $(X, g_+)$ satisfies $Ric_{g_+}=-dg_+$, $(X, g_+)$ is called a conformally compact Einstein manifold.
For a conformally compact Einstein metric and any defining function $r$, $|dr|^2_{g_c}=1$ on $M$.
Hence conformally compact Einstein manifolds are asymptotic hyperbolic Einstein manifolds.
From now on let $(X^{d+1}, g_+)$ be a conformally compact Einstein manifold.
For a conformally compact Einstein manifold $(X, g_+)$, as indicated in the introduction,
there is a one-to-one correspondence between $\overline{g}\in [g_{confinf}]$ and special defining functions defined near $\partial X$,
and the identification (\ref{identification}) of $M\times [0, \epsilon)$ with a neighborhood of $M$ in $\overline{X}$ so that
$g_+=r^{-2}(dr^2+\widetilde{g})$,
where an expansion of $\widetilde{g}$ is given by (\ref{oddexpansion}) and (\ref{evenexpansion}) for $d$ odd and even respectively.
The determined coefficients in these expansions can be calculated and are combinations of $\overline{g}$ and
its derivatives, see \cite{FeffermanGraham, Graham09}. For example, for $d\geq 3$ one has
\begin{equation}\label{g2}
g^{(2)}=-\overline{P}=-\frac{1}{d-2}[\overline{Ric}-\frac{\overline{R}}{2(d-1)}\overline{g}],
\end{equation}
here $\overline{P}$ is the Schouten tensor of $\overline{g}$. It is well-known that
\begin{equation}\label{WeylSchouten}
\overline{R}_{ABCD}=\overline{W}_{ABCD}+\overline{P}_{AC}\overline{g}_{BD}+\overline{P}_{BD}\overline{g}_{AC}
-\overline{P}_{AD}\overline{g}_{BC}-\overline{P}_{BC}\overline{g}_{AD},
\end{equation}
here $\overline{W}$ is the Weyl tensor of $\overline{g}$. For $d\geq 5$,
\begin{equation}\label{g4}
g^{(4)}_{CD}=\frac{1}{4(4-d)}\overline{B}_{CD}+\frac{1}{4}\overline{g}^{EF}\overline{P}_{CE}\overline{P}_{DF},
\end{equation}
where
\begin{equation}\label{Bach}
\overline{B}_{CD}=\overline{\triangle} \overline{P}_{CD}-\overline{\nabla}^E\overline{\nabla}_DP_{CE}+\overline{P}^{EF}\overline{W}_{CEDF}
\end{equation}
is the Bach tensor of $\overline{g}$.

Let $(X^{d+1},g_+)$ be a conformally compact Einstein manifold and $r$ a special defining function such that $g_c|_M=\overline{g}$.
Now we consider a minimal surface $Y^{n+1}$, $0\leq n\leq d-1$, of $(X, g_+)$ with the smooth boundary $\Sigma^n\hookrightarrow M^d$.
Near a point of $\Sigma$, let $(x^i, y^\alpha)$ be a local coordinate chart as indicated in (\ref{coordinatexy}),
which is extended into a neighborhood of the point in $\overline{X}$ via the identification (\ref{identification}).
The second fundamental form of $\Sigma\hookrightarrow (M, \overline{g})$ is
\begin{equation}\label{secondfundamentalform}
A_{ij}^\alpha=\overline{g}^{\alpha\beta}\overline{g}(\overline{\nabla}_i\partial_j,\partial_\beta)
=-\frac{1}{2}\overline{g}^{\alpha\beta}\partial_\beta \overline{g}_{ij}=\overline{\Gamma}_{ij}^\alpha.
\end{equation}

We consider minimal surfaces $Y$ in $(X, g_+)$ which may be written as a graph $\{y^\alpha=u^\alpha(x,r)\}$.
That is near the boundary, $Y=(x^i,u^\alpha(x,r),r)$. Let
\begin{equation}\label{h}
h=g_+|_Y, \quad \overline{h}=r^2h=g_c|_Y=(dr^2+\widetilde{g})|_Y,
\end{equation}
where $\widetilde{g}$ is given by (\ref{gr}).
Near the boundary, the tangent space of $Y$ is spanned by
$$Y_r=\partial_r+u_r^\gamma \partial_\gamma,$$
and
$$Y_i=\partial_i+u_i^\gamma\partial_\gamma, \quad i=1,2,\cdots,n.$$
Then
$$\overline{h}_{rr}=\overline{h}(Y_r,Y_r)=1+\widetilde{g}_{\alpha\beta}u_r^\alpha u_r^\beta,$$
$$\overline{h}_{ir}=\overline{h}(Y_i,Y_r)=\widetilde{g}_{i\alpha}u_r^\alpha+\widetilde{g}_{\alpha\beta}u_i^\alpha u_r^\beta,$$
and
$$\overline{h}_{ij}=\overline{h}(Y_i,Y_j)=\widetilde{g}_{ij}+\widetilde{g}_{i\alpha}u_j^\alpha
+\widetilde{g}_{j\alpha}u_i^\alpha+\widetilde{g}_{\alpha\beta}u_i^\alpha u_j^\beta.$$
Graham and Witten \cite{GrahamWitten} showed that the graphic minimal surface equation of $Y$ is
$$M(u)=0,$$
where
\begin{eqnarray}\label{minimal}
M(u)_\gamma&=&[r\partial_r-(n+1)+\frac{1}{2}r\partial_r\overline{L}][\overline{h}^{rr}\widetilde{g}_{\beta\gamma}u_r^\beta
+\overline{h}^{ir}(\widetilde{g}_{i\gamma}+\widetilde{g}_{\beta\gamma}u_i^\beta)]
\nonumber
\\&&+r[\partial_j+\frac{1}{2}\partial_j\overline{L}][\overline{h}^{rj}\widetilde{g}_{\beta\gamma}u_r^\beta+
\overline{h}^{ij}(\widetilde{g}_{i\gamma}+\widetilde{g}_{\beta\gamma}u_i^\beta)]\nonumber
\\&&-\frac{1}{2}r\overline{h}^{ij}[\partial_\gamma \widetilde{g}_{ij}+2\partial_\gamma \widetilde{g}_{i\alpha}u_j^\alpha+\partial_\gamma \widetilde{g}_{\alpha\beta}u_i^\alpha u_j^\beta]
\\&&-r\overline{h}^{ir}[\partial_\gamma \widetilde{g}_{i\alpha}u_r^\alpha+\partial_\gamma \widetilde{g}_{\alpha\beta}u_i^\alpha u_r^\beta]\nonumber
\\&&-\frac{1}{2}r\overline{h}^{rr}\partial_\gamma \widetilde{g}_{\alpha\beta}u_r^\alpha u_r^\beta,\nonumber
\end{eqnarray}
and
\begin{equation}\label{Lbar}
\overline{L}=\log \det \overline{h}.
\end{equation}
By the choice of the local coordinates $(x^i,y^\alpha)$ on $M$, we have $u(x,0)=0$.
Graham and Witten \cite{GrahamWitten} found that at $r=0$,
\begin{equation}\label{ur}
u_r=0,
\end{equation}
which means that $Y$ intersects with $M$ orthogonally; and for $n\geq 1$
\begin{equation}\label{urr}
u_{rr}^\alpha=-\frac{1}{2n}g^{ij}\overline{g}^{\alpha\beta}\partial_\beta \overline{g}_{ij}=\frac{1}{n}H^\alpha,
\end{equation}
here $H$ is the mean curvature of $\Sigma\hookrightarrow (M, \overline{g})$.
The expansion of the height functions $u$ in $r$, up to a critical order, takes the form of (\ref{uexpansionoddn}) for $n$ odd
and of (\ref{uexpansionevenn}) for $n$ even. Namely, for $n$ odd
$$u=u^{(2)}r^2+(\mbox{even powers})+u^{(n+1)}r^{n+1}+u^{(n+2)}r^{n+2}+\cdots,$$
and for $n$ even
$$u=u^{(2)}r^2+(\mbox{even powers})+u^{(n)}r^{n}+w_nr^{n+2}\log r+u^{(n+2)}r^{n+2}+\cdots,$$
where $u^{(k)}, k<n+2,$ and $w_n$ are locally determined functions of $\Sigma$.

\subsection{Coefficient $u^{(4)}$ in (\ref{uexpansionoddn}) and (\ref{uexpansionevenn}) for $n\geq 3$.}
For $n\geq 3$, let
\begin{equation}\label{vw}
u=vr^2+wr^4+O(r^5),
\end{equation}
where, it follows from (\ref{urr}),
\begin{equation}\label{vH}
v^\alpha=\frac{1}{2n}H^\alpha.
\end{equation}
It follows from (\ref{Lbar}) that
$$\partial_r\overline{L}=O(r). %\quad \quad \partial_iL=\partial_i\log\det(g|_{r=0})+O(r^2).
$$
Note that
$$\overline{h}_{ir}=O(r^3), \quad \overline{h}^{ir}=O(r^3),$$
$$\widetilde{g}_{i\alpha}=O(r^2),  \quad u_r^\beta=O(r), \quad u_i^\beta=O(r^2),$$
and the fact that the coefficient $w$ in (\ref{vw}) can be derived by tracing the coefficient of $r^3$ in (\ref{minimal}).
First, we have
\begin{eqnarray}
M(u)_\gamma&=&[r\partial_r-(n+1)](\overline{h}^{rr}\widetilde{g}_{\beta\gamma}u_r^\beta)
+\frac{1}{2}r\partial_r\overline{L}\overline{h}^{rr}\widetilde{g}_{\beta\gamma}u_r^\beta \label{a}
\\&&+r\partial_j[\overline{h}^{ij}(\widetilde{g}_{i\gamma}+\widetilde{g}_{\beta\gamma}u_i^\beta)]
+\frac{1}{2}r\partial_j\overline{L}\overline{h}^{ij}(\widetilde{g}_{i\gamma}+\widetilde{g}_{\beta\gamma}u_i^\beta)\label{b}
\\&&-\frac{1}{2}r\overline{h}^{ij}\partial_\gamma \widetilde{g}_{ij}-r\overline{h}^{ij}\partial_\gamma \widetilde{g}_{i\alpha}u_j^\alpha  \label{c}
-\frac{1}{2}r\overline{h}^{rr}\partial_\gamma \widetilde{g}_{\alpha\beta}u_r^\alpha u_r^\beta\label{d}
\\&&+o(r^3)\nonumber.
\end{eqnarray}

\begin{lemma} For the right hand side of (\ref{a}), we have
\begin{eqnarray}\label{aA}
\lefteqn{[r\partial_r-(n+1)](\overline{h}^{rr}\widetilde{g}_{\beta\gamma}u_r^\beta)
+\frac{1}{2}r\partial_r\overline{L}\overline{h}^{rr}\widetilde{g}_{\beta\gamma}u_r^\beta}\nonumber
\\&=&-2nr\overline{g}_{\beta\gamma}v^\beta-(n-2)(2g^{(2)}_{\beta\gamma}v^\beta
+2\partial_\alpha \overline{g}_{\beta \gamma}v^\alpha v^\beta+4\overline{g}_{\beta\gamma}w^\beta)r^3
\\&&+r^3[-8|v|^2\overline{g}_{\beta\gamma}v^\beta+2g^{ij}g^{(2)}_{ij}\overline{g}_{\beta\gamma}v^\beta]+o(r^3).\nonumber
\end{eqnarray}
\end{lemma}
Proof. Note that
$$\overline{h}^{rr}=1-4|v|^2r^2+o(r^2).$$
Hence
\begin{eqnarray*}
\overline{h}^{rr}\widetilde{g}_{\beta\gamma}u_r^\beta&=&(1-4|v|^2r^2)(\overline{g}_{\beta\gamma}
+r^2g^{(2)}_{\beta\gamma}+r^2\partial_\alpha \overline{g}_{\beta \gamma}v^\alpha)(2rv^\beta+4r^3w^\beta)+o(r^3)
\\&=&2r\overline{g}_{\beta\gamma}v^\beta+(-8|v|^2\overline{g}_{\beta\gamma}v^\beta +2g^{(2)}_{\beta\gamma}v^\beta
+2\partial_\alpha \overline{g}_{\beta \gamma}v^\alpha v^\beta+4\overline{g}_{\beta\gamma}w^\beta)r^3+o(r^3).
\end{eqnarray*}
Therefore,
\begin{eqnarray}\label{lemma1-1}
\lefteqn{[r\partial_r-(n+1)](\overline{h}^{rr}\widetilde{g}_{\beta\gamma}u_r^\beta)
=-2nr\overline{g}_{\beta\gamma}v^\beta} \nonumber
\\&&-(n-2)(-8|v|^2\overline{g}_{\beta\gamma}v^\beta +2g^{(2)}_{\beta\gamma}v^\beta
+2\partial_\alpha \overline{g}_{\beta \gamma}v^\alpha v^\beta+4\overline{g}_{\beta\gamma}w^\beta)r^3+o(r^3).
\end{eqnarray}

Note that
\begin{eqnarray*}
\partial_r\overline{L}&=&\overline{h}^{rr}\partial_r\overline{h}_{rr}+\overline{h}^{ij}\partial_r\overline{h}_{ij}+o(r)
\\&=&8r\overline{g}_{\alpha\beta}v^\alpha v^\beta+g^{ij}(2rg^{(2)}_{ij}+2r\partial_\alpha \overline{g}_{ij} v^\alpha)+o(r)
\\&=&r[(8-8n)|v|^2+2g^{ij}g^{(2)}_{ij}]+o(r),
\end{eqnarray*}
where in the last equality we used (\ref{secondfundamentalform}) and (\ref{vH}).
Hence
\begin{equation}\label{lemma1-2}
\frac{1}{2}r\partial_r\overline{L}\overline{h}^{rr}\widetilde{g}_{\beta\gamma}u_r^\beta=
r^3[(8-8n)|v|^2+2g^{ij}g^{(2)}_{ij}]\overline{g}_{\beta\gamma}v^\beta+o(r^3).
\end{equation}
(\ref{aA}) then follows from (\ref{lemma1-1}) and (\ref{lemma1-2}).  $\hfill \Box$

\begin{lemma} For (\ref{b}), we have
\begin{eqnarray}\label{bA}
\lefteqn{
r\partial_j[\overline{h}^{ij}(\widetilde{g}_{i\gamma}+\widetilde{g}_{\beta\gamma}u_i^\beta)]
+\frac{1}{2}r\partial_j\overline{L}\overline{h}^{ij}(\widetilde{g}_{i\gamma}+\widetilde{g}_{\beta\gamma}u_i^\beta)}\nonumber
\\&=&r^3g^{ij}[\partial_jg^{(2)}_{i\gamma}+\partial_j\partial_\beta \overline{g}_{i\gamma}v^\beta
+\partial_\beta \overline{g}_{i\gamma}v_j^\beta+\partial_j\overline{g}_{\beta\gamma}v_i^\beta+\overline{g}_{\beta \gamma}\partial_jv_i^\beta]
\\&&-r^3g^{kl}\Gamma_{kl}^i(g^{(2)}_{i\gamma}+\partial_\beta\overline{ g}_{i\gamma}v^\beta+\overline{g}_{\beta\gamma}v_i^\beta)+o(r^3).\nonumber
\end{eqnarray}
\end{lemma}
Proof. It is easy to see that
\begin{eqnarray*}
r\partial_j[\overline{h}^{ij}(\widetilde{g}_{i\gamma}+\widetilde{g}_{\beta\gamma}u_i^\beta)]
&=&-r^3g^{ik}g^{jl}\partial_jg_{kl}(g^{(2)}_{i\gamma}+\partial_\beta \overline{g}_{i\gamma}v^\beta+\overline{g}_{\beta\gamma}v_i^\beta)
\\&&+rg^{ij}\partial_j\widetilde{g}_{i\gamma}+r^3g^{ij}(\partial_j\overline{g}_{\beta \gamma}v_i^\beta+\overline{g}_{\beta \gamma}\partial_jv_i^\beta)+o(r^3).
\end{eqnarray*}
It then follows from
\begin{eqnarray*}
\partial_j\widetilde{g}_{i\gamma}&=&r^2(\partial_jg^{(2)}_{i\gamma}+\partial_j\partial_\beta \overline{g}_{i\gamma}v^\beta
+\partial_\beta \overline{g}_{i\gamma}v_j^\beta)+o(r^2)
\end{eqnarray*}
that
\begin{eqnarray}\label{lemma2-1}
\lefteqn{r\partial_j[\overline{h}^{ij}(\widetilde{g}_{i\gamma}+\widetilde{g}_{\beta\gamma}u_i^\beta)]
=-r^3g^{ik}g^{jl}\partial_jg_{kl}(g^{(2)}_{i\gamma}+\partial_\beta \overline{g}_{i\gamma}v^\beta+\overline{g}_{\beta\gamma}v_i^\beta)}\nonumber
\\&&+r^3g^{ij}[\partial_jg^{(2)}_{i\gamma}+\partial_j\partial_\beta \overline{g}_{i\gamma}v^\beta
+\partial_\beta \overline{g}_{i\gamma}v_j^\beta+\partial_j\overline{g}_{\beta\gamma}v_i^\beta+\overline{g}_{\beta \gamma}\partial_jv_i^\beta]+o(r^3).
\end{eqnarray}

It is easy to see that
\begin{eqnarray}\label{lemma2-2}
\lefteqn{\frac{1}{2}r\partial_j\overline{L}\overline{h}^{ij}(\widetilde{g}_{i\gamma}+\widetilde{g}_{\beta\gamma}u_i^\beta)}\nonumber
\\&=&\frac{1}{2}r^3g^{ij}g^{kl}\partial_jg_{kl}(g^{(2)}_{i\gamma}+\partial_\beta \overline{g}_{i\gamma}v^\beta+\overline{g}_{\beta\gamma}v_i^\beta)+o(r^3).
\end{eqnarray}

(\ref{bA}) then follows from (\ref{lemma2-1}) and (\ref{lemma2-2}). $\hfill \Box$

\begin{lemma} For (\ref{c}), we have
\begin{eqnarray}\label{cA}
\lefteqn{-\frac{1}{2}r\overline{h}^{ij}\partial_\gamma \widetilde{g}_{ij}-r\overline{h}^{ij}\partial_\gamma \widetilde{g}_{i\alpha}u_j^\alpha
-\frac{1}{2}r\overline{h}^{rr}\partial_\gamma \widetilde{g}_{\alpha\beta}u_r^\alpha u_r^\beta}\nonumber
\\&&=2nr\overline{g}_{\beta\gamma}v^\beta+[-\frac{1}{2}g^{ij}\partial_\gamma g^{(2)}_{ij}-\frac{1}{2}g^{ij}\partial_\gamma\partial_\beta \overline{g}_{ij}v^\beta
-g^{ik}g^{jl}g^{(2)}_{kl}A_{ij}^\beta \overline{g}_{\beta\gamma}]r^3
\\&&+2r^3g^{ik}g^{jl}\langle A_{kl},v\rangle A_{ij}^\beta \overline{g}_{\beta\gamma}
+r^3[-g^{ij}\partial_\gamma \overline{g}_{i\alpha}v_j^\alpha -2\partial_\gamma \overline{g}_{\alpha\beta}v^\alpha v^\beta]+o(r^3).\nonumber
\end{eqnarray}
\end{lemma}
Proof. It is clear that we have
\begin{eqnarray*}
\overline{h}^{ij}&=&g^{ij}-g^{ik}g^{jl}(g^{(2)}_{kl}+\partial_\alpha \overline{g}_{kl} v^\alpha)r^2+o(r^2),
\end{eqnarray*}
and
$$\partial_\gamma \widetilde{g}_{ij}=\partial_\gamma \overline{g}_{ij}
+r^2(\partial_\gamma g^{(2)}_{ij}+\partial_\gamma\partial_\beta \overline{g}_{ij} v^\beta)+o(r^2).$$
Hence by using (\ref{secondfundamentalform}), we get
\begin{eqnarray}\label{Lemma3-1}
\lefteqn{-\frac{1}{2}r\overline{h}^{ij}\partial_\gamma \widetilde{g}_{ij}=2nr\overline{g}_{\beta\gamma}v^\beta}\nonumber
\\&&+[-\frac{1}{2}g^{ij}\partial_\gamma g^{(2)}_{ij}-\frac{1}{2}g^{ij}\partial_\gamma\partial_\beta \overline{g}_{ij}v^\beta
-g^{ik}g^{jl}g^{(2)}_{kl}A_{ij}^\beta \overline{g}_{\beta\gamma}]r^3
\\&&+2r^3g^{ik}g^{jl}\langle A_{kl},v\rangle A_{ij}^\beta \overline{g}_{\beta\gamma}+o(r^3).\nonumber
\end{eqnarray}

On the other hand, it's easy to see that
\begin{eqnarray}\label{Lemma3-2}
\lefteqn{-r\overline{h}^{ij}\partial_\gamma \widetilde{g}_{i\alpha}u_j^\alpha
-\frac{1}{2}r\overline{h}^{rr}\partial_\gamma \widetilde{g}_{\alpha\beta}u_r^\alpha u_r^\beta } \nonumber
\\&=&r^3[-g^{ij}\partial_\gamma \overline{g}_{i\alpha}v_j^\alpha
-2\partial_\gamma \overline{g}_{\alpha\beta}v^\alpha v^\beta]+o(r^3).
\end{eqnarray}

(\ref{cA}) then follows from (\ref{Lemma3-1}) and (\ref{Lemma3-2}). $\hfill \Box$

Putting (\ref{aA}) (\ref{bA}) and (\ref{cA}) together, we get the following
\begin{corollary}
For $n\geq 3$, we have
\begin{eqnarray}
M(u)_\gamma&=&-(n-2)(2g^{(2)}_{\beta\gamma}v^\beta
+2\partial_\alpha \overline{g}_{\beta \gamma}v^\alpha v^\beta+4\overline{g}_{\beta\gamma}w^\beta)r^3 \nonumber
\\&&+r^3[-8|v|^2\overline{g}_{\beta\gamma}v^\beta+2g^{ij}g^{(2)}_{ij}\overline{g}_{\beta\gamma}v^\beta] \nonumber
\\&&+r^3g^{ij}[\partial_jg^{(2)}_{i\gamma}+\partial_j\partial_\beta \overline{g}_{i\gamma}v^\beta
+\partial_\beta \overline{g}_{i\gamma}v_j^\beta+\partial_j\overline{g}_{\beta\gamma}v_i^\beta+\overline{g}_{\beta \gamma}\partial_jv_i^\beta] \nonumber
\\&&-r^3g^{kl}\Gamma_{kl}^i(g^{(2)}_{i\gamma}+\partial_\beta\overline{ g}_{i\gamma}v^\beta+\overline{g}_{\beta\gamma}v_i^\beta)\label{abc}
\\&&+[-\frac{1}{2}g^{ij}\partial_\gamma g^{(2)}_{ij}-\frac{1}{2}g^{ij}\partial_\gamma\partial_\beta \overline{g}_{ij}v^\beta
-g^{ik}g^{jl}g^{(2)}_{kl}A_{ij}^\beta \overline{g}_{\beta\gamma}]r^3 \nonumber
\\&&+2r^3g^{ik}g^{jl}\langle A_{kl},v\rangle A_{ij}^\beta \overline{g}_{\beta\gamma}
+r^3[-g^{ij}\partial_\gamma \overline{g}_{i\alpha}v_j^\alpha -2\partial_\gamma \overline{g}_{\alpha\beta}v^\alpha v^\beta]+o(r^3).\nonumber
\end{eqnarray}
\end{corollary}

Note that
$$w=u^{(4)},$$
so we have
\begin{proposition}
Let $Y^{n+1}$, $n\geq 3$, be a minimal surface in $(X^{d+1}, g_+)$. Then
\begin{eqnarray}\label{w=}
4(n-2)\overline{g}_{\beta\gamma}u^{(4)\beta}&=&-2(n-2)g^{(2)}_{\beta\gamma}v^\beta-8|v|^2\overline{g}_{\beta\gamma}v^\beta
+2g^{ij}g^{(2)}_{ij}\overline{g}_{\beta\gamma}v^\beta \nonumber
\\&&-g^{ik}g^{jl}g^{(2)}_{kl}A_{ij}^\beta \overline{g}_{\beta\gamma}+2g^{ik}g^{jl}\langle A_{kl},v\rangle A_{ij}^\beta \overline{g}_{\beta\gamma}
\\&&+Q_\gamma,\nonumber
\end{eqnarray}
where
\begin{eqnarray}
Q_\gamma&=&g^{ij}[\partial_jg^{(2)}_{i\gamma}+\partial_j\partial_\beta \overline{g}_{i\gamma}v^\beta+\partial_\beta \overline{g}_{i\gamma}v_j^\beta
+\partial_j\overline{g}_{\beta\gamma}v_i^\beta-\partial_\gamma \overline{g}_{i\alpha}v_j^\alpha+\overline{g}_{\beta \gamma}\partial_jv_i^\beta] \nonumber
\\&&-g^{kl}\Gamma_{kl}^i(g^{(2)}_{i\gamma}+\partial_\beta \overline{g}_{i\gamma}v^\beta+\overline{g}_{\beta\gamma}v_i^\beta) \nonumber
\\&&-\frac{1}{2}g^{ij}\partial_\gamma g^{(2)}_{ij}-\frac{1}{2}g^{ij}\partial_\gamma\partial_\beta \overline{g}_{ij}v^\beta \label{Q}
\\&&-2(n-2)\partial_\alpha \overline{g}_{\beta \gamma}v^\alpha v^\beta-2\partial_\gamma \overline{g}_{\alpha\beta}v^\alpha v^\beta.\nonumber
\end{eqnarray}
\end{proposition}

\section{Proof of Proposition \ref{u4prop}}

Assume $n\geq 3$.
In this section, we will reformulate $Q_\gamma$ to get the expression of $u^{(4)}$, as given by (\ref{u4intro}).
Let
$$Q^\alpha=\overline{g}^{\alpha\gamma}Q_\gamma.$$
It follows from (\ref{Q}) that
\begin{eqnarray}\label{w2=}
4(n-2)u^{(4)\alpha}&=&-2(n-2)\overline{g}^{\alpha\gamma}g^{(2)}_{\beta\gamma}v^\beta
-8|v|^2v^\alpha+2g^{ij}g^{(2)}_{ij}v^\alpha \nonumber
\\&&-g^{ik}g^{jl}g^{(2)}_{kl}A_{ij}^\alpha+2g^{ik}g^{jl}\langle A_{kl},v\rangle A_{ij}^\alpha+Q^\alpha,
\end{eqnarray}
where
\begin{eqnarray*}
Q^\alpha&=&I^\alpha+II^\alpha,
\end{eqnarray*}
and
\begin{eqnarray}
I^\alpha&=&g^{ij}[\partial_jv_i^\alpha-\frac{1}{2}\overline{g}^{\alpha\gamma}\partial_\gamma\partial_\beta \overline{g}_{ij}v^\beta
+\overline{g}^{\alpha\gamma}\partial_j\partial_\beta \overline{g}_{i\gamma}v^\beta]
+2g^{ij}\overline{\Gamma}^\alpha_{i\beta}v_j^\beta\nonumber
\\&&-2(n-2)\overline{g}^{\alpha\gamma}\partial_\beta \overline{g}_{\eta\gamma}v^\beta v^\eta-2\overline{g}^{\alpha\gamma}\partial_\gamma \overline{g}_{\beta\eta}v^\beta v^\eta \label{I}
\\&&-\overline{g}^{\alpha\gamma}g^{kl}\Gamma_{kl}^i(\partial_\beta \overline{g}_{i\gamma}v^\beta+\overline{g}_{\beta\gamma}v_i^\beta),\nonumber
\end{eqnarray}
$$II^\alpha=g^{ij}\overline{g}^{\alpha\gamma}[-\frac{1}{2}\partial_\gamma g^{(2)}_{ij}+\partial_jg^{(2)}_{i\gamma}-\Gamma_{ij}^kg^{(2)}_{k\gamma}].$$

We first compute $II^\alpha$.

\begin{lemma}\label{II} We have
\begin{eqnarray}\label{IIformula}
II^\alpha&=&-\frac{1}{2}\overline{g}^{\alpha\gamma}g^{ij}(\overline{\nabla}_\gamma g^{(2)}_{ij}-2\overline{\nabla}_ig^{(2)}_{j\gamma})
+\overline{g}^{\alpha\gamma}H^\beta g^{(2)}_{\beta\gamma}.
\end{eqnarray}
\end{lemma}
Proof. Let $D$ range from $1$ to $d$ so that for $\partial_D=\partial_i$ for $D=i\leq n$ and $\partial_D=\partial_\alpha$ for $D=\alpha\geq n+1$.
We have
\begin{eqnarray}
\lefteqn{\partial_\gamma g^{(2)}_{ij}-\partial_jg^{(2)}_{i\gamma}-\partial_ig^{(2)}_{j\gamma}}\nonumber
\\&=&(\overline{\nabla}_\gamma g^{(2)}_{ij}+\overline{\Gamma}_{\gamma i}^Dg^{(2)}_{Dj}+\overline{\Gamma}_{\gamma j}^Dg^{(2)}_{iD})
-(\overline{\nabla}_jg^{(2)}_{i\gamma}+\overline{\Gamma}_{ji}^Dg^{(2)}_{D\gamma}+\overline{\Gamma}_{j\gamma}^Dg^{(2)}_{iD})\nonumber
\\&&-(\overline{\nabla}_ig^{(2)}_{j\gamma}+\overline{\Gamma}_{ij}^Dg^{(2)}_{D\gamma}+\overline{\Gamma}_{i\gamma}^Dg^{(2)}_{jD})\nonumber
\\&=&\overline{\nabla}_\gamma g^{(2)}_{ij}-\overline{\nabla}_jg^{(2)}_{i\gamma}-\overline{\nabla}_ig^{(2)}_{j\gamma}-2\overline{\Gamma}_{ij}^Dg^{(2)}_{D\gamma}.\label{II}
\end{eqnarray}
Hence by using (\ref{secondfundamentalform}), we get
\begin{eqnarray*}
II^\alpha&=&g^{ij}\overline{g}^{\alpha\gamma}[-\frac{1}{2}\partial_\gamma g^{(2)}_{ij}+\partial_jg^{(2)}_{i\gamma}-\Gamma_{ij}^kg^{(2)}_{k\gamma}]
\\&=&-\frac{1}{2}\overline{g}^{\alpha\gamma}g^{ij}(\overline{\nabla}_\gamma g^{(2)}_{ij}-2\overline{\nabla}_ig^{(2)}_{j\gamma})+\overline{g}^{\alpha\gamma}g^{ij}\overline{\Gamma}_{ij}^Dg^{(2)}_{D\gamma}
-\overline{g}^{\alpha\gamma}g^{ij}\Gamma_{ij}^kg^{(2)}_{k\gamma}
\\&=&-\frac{1}{2}\overline{g}^{\alpha\gamma}g^{ij}(\overline{\nabla}_\gamma g^{(2)}_{ij}-2\overline{\nabla}_ig^{(2)}_{j\gamma})
+\overline{g}^{\alpha\gamma}H^\beta g^{(2)}_{\beta\gamma}.
\end{eqnarray*}
$\hfill \Box$

We now deal with $I^\alpha$.

\begin{lemma}
We have
\begin{eqnarray}
\partial_jv_i^\alpha&=&\nabla_j^\perp\nabla_i^\perp v^\alpha-g^{kl}\langle A_{ik},v\rangle A_{jl}^\alpha
-v_i^\beta \overline{\Gamma}_{j\beta}^\alpha-v_j^\beta \overline{\Gamma}_{i\beta}^\alpha+v_k^\alpha\Gamma_{ij}^k \nonumber
\\&&-v^\beta(\partial_j\overline{\Gamma}_{i\beta}^\alpha+\overline{\Gamma}_{i\beta}^\gamma\overline{\Gamma}_{j\gamma}^\alpha
+\overline{\Gamma}_{i\beta}^k\overline{\Gamma}_{jk}^\alpha-\Gamma_{ij}^k\overline{\Gamma}_{k\beta}^\alpha). \label{Ia}
\end{eqnarray}
\end{lemma}
Proof. Let $A, B, C, D=\{i,\alpha\}=1,2,\cdots,d$.
Let $X^\perp$ denote the normal part of $X$ and $\nabla^\perp$ the covariant differentiation with respect to the normal connection.
We have
\begin{eqnarray*}
\overline{\nabla}_iv&=&\partial_iv^\alpha\partial_\alpha+v^\beta\overline{\Gamma}_{i\beta}^A\partial_A
\\&=&(v_i^\alpha+v^\beta\overline{\Gamma}_{i\beta}^\alpha)\partial_\alpha+v^\beta\overline{\Gamma}_{i\beta}^k\partial_k.
\end{eqnarray*}
Then
\begin{eqnarray*}
(\overline{\nabla}_j\overline{\nabla}_iv)^\perp&=&[\overline{\nabla}_j(\overline{\nabla}_iv)]^\perp-\Gamma_{ij}^k\nabla_k^\perp v
\\&=&\partial_j(v_i^\alpha+v^\beta\overline{\Gamma}_{i\beta}^\alpha)\partial_\alpha
+v_i^\beta \overline{\Gamma}_{j\beta}^\alpha\partial_\alpha+v^\beta \overline{\Gamma}_{i\beta}^D\overline{\Gamma}_{Dj}^\alpha\partial_\alpha
-\Gamma_{ij}^k\nabla_k^\perp v
\\&=&(\partial_jv_i^\alpha+v_i^\beta \overline{\Gamma}_{j\beta}^\alpha+v_j^\beta \overline{\Gamma}_{i\beta}^\alpha-v_k^\alpha\Gamma_{ij}^k)\partial_\alpha
\\&&+v^\beta(\partial_j\overline{\Gamma}_{i\beta}^\alpha+\overline{\Gamma}_{i\beta}^D\overline{\Gamma}_{Dj}^\alpha
-\Gamma_{ij}^k\overline{\Gamma}_{k\beta}^\alpha)\partial_\alpha.
\end{eqnarray*}
%note that here $\overline{\nabla}_j\overline{\nabla}_iv$ is understood as $\overline{\nabla}_j(\overline{\nabla}_iv)-\Gamma_{ij}^k\overline{\nabla}_kv$, i.e.
%the differentiation with respect to the bundle $i^*TM$ with $i:N\rightarrow M$.
On the other hand,
\begin{eqnarray*}
\nabla_j^\perp\nabla_i^\perp v&=&[\overline{\nabla}_j(\nabla_i^\perp v)]^\perp-\Gamma_{ij}^k\nabla_k^\perp v
\\&=&[\overline{\nabla}_j(\overline{\nabla}_iv)]^\perp+g^{kl}\langle A_{ik},v\rangle A_{jl}-\Gamma_{ij}^k\nabla_k^\perp v
\\&=&(\overline{\nabla}_j\overline{\nabla}_iv)^\perp+g^{kl}\langle A_{ik},v\rangle A_{jl}.
\end{eqnarray*}
Therefore, (\ref{Ia}) follows from
\begin{eqnarray*}
\nabla_j^\perp\nabla_i^\perp v&=&(\partial_jv_i^\alpha+v_i^\beta \overline{\Gamma}_{j\beta}^\alpha+v_j^\beta \overline{\Gamma}_{i\beta}^\alpha-v_k^\alpha\Gamma_{ij}^k)\partial_\alpha
\\&&+v^\beta(\partial_j\overline{\Gamma}_{i\beta}^\alpha+\overline{\Gamma}_{i\beta}^D\overline{\Gamma}_{Dj}^\alpha
-\Gamma_{ij}^k\overline{\Gamma}_{k\beta}^\alpha)\partial_\alpha
\\&&+g^{kl}\langle A_{ik},v\rangle A_{jl}.
\end{eqnarray*}
$\hfill \Box$

\begin{lemma} We have
\begin{eqnarray}
\lefteqn{g^{ij}[-\frac{1}{2}\overline{g}^{\alpha\gamma}\partial_\gamma\partial_\beta \overline{g}_{ij}v^\beta
+\overline{g}^{\alpha\gamma}\partial_j\partial_\beta \overline{g}_{i\gamma}v^\beta]}\nonumber
\\&=&g^{ij}v^\beta\partial_\beta\overline{\Gamma}_{ij}^\alpha+g^{ij}\Gamma_{ij}^l\overline{g}^{\alpha\xi}v^\beta\partial_\beta \overline{g}_{l\xi}
+2n\overline{g}^{\alpha\xi} v^\beta v^\eta  \partial_\beta \overline{g}_{\xi\eta}.\label{Ib}
\end{eqnarray}
\end{lemma}
Proof. Note that
$$\overline{\Gamma}_{ij}^\alpha=\frac{1}{2}\overline{g}^{\alpha\gamma}(\partial_i\overline{g}_{j\gamma}+\partial_j\overline{g}_{i\gamma}-\partial_\gamma \overline{g}_{ij})
+\frac{1}{2}\overline{g}^{\alpha k}(\partial_i\overline{g}_{jk}+\partial_j\overline{g}_{ik}-\partial_k \overline{g}_{ij}),$$
hence
\begin{eqnarray*}
\lefteqn{g^{ij}[-\frac{1}{2}\overline{g}^{\alpha\gamma}\partial_\gamma\partial_\beta \overline{g}_{ij}v^\beta
+\overline{g}^{\alpha\gamma}\partial_j\partial_\beta \overline{g}_{i\gamma}v^\beta]}
\\&=&g^{ij}v^\beta(\partial_\beta[\overline{g}^{\alpha\gamma}(-\frac{1}{2}\partial_\gamma \overline{g}_{ij}+\partial_j\overline{g}_{i\gamma})]
-\partial_\beta \overline{g}^{\alpha\gamma}(-\frac{1}{2}\partial_\gamma \overline{g}_{ij}+\partial_j\overline{g}_{i\gamma}))
\\&=&g^{ij}v^\beta(\partial_\beta[\overline{\Gamma}_{ij}^\alpha-\frac{1}{2}\overline{g}^{\alpha k}(\partial_i\overline{g}_{jk}+\partial_j\overline{g}_{ik}-\partial_k \overline{g}_{ij})]
-\partial_\beta \overline{g}^{\alpha\gamma}(-\frac{1}{2}\partial_\gamma \overline{g}_{ij}+\partial_j\overline{g}_{i\gamma})).
\end{eqnarray*}
Note that on $\Sigma^n$
$$\overline{g}^{\alpha k}=0,$$
hence we have
\begin{eqnarray*}
\lefteqn{g^{ij}[-\frac{1}{2}\overline{g}^{\alpha\gamma}\partial_\gamma\partial_\beta \overline{g}_{ij}v^\beta
+\overline{g}^{\alpha\gamma}\partial_j\partial_\beta \overline{g}_{i\gamma}v^\beta]}
\\&=&g^{ij}v^\beta\partial_\beta\overline{\Gamma}_{ij}^\alpha
-\frac{1}{2}g^{ij}v^\beta\partial_\beta \overline{g}^{\alpha k}(\partial_i\overline{g}_{jk}+\partial_j\overline{g}_{ik}-\partial_k \overline{g}_{ij})
\\&&-g^{ij}v^\beta\partial_\beta \overline{g}^{\alpha\gamma}(-\frac{1}{2}\partial_\gamma \overline{g}_{ij}+\partial_j\overline{g}_{i\gamma})
\\&=&g^{ij}v^\beta\partial_\beta\overline{\Gamma}_{ij}^\alpha
+\frac{1}{2}g^{ij}v^\beta \overline{g}^{\alpha\xi}g^{kl}\partial_\beta \overline{g}_{l\xi }(\partial_ig_{jk}+\partial_jg_{ik}-\partial_k g_{ij})
\\&&+g^{ij}v^\beta \overline{g}^{\alpha\xi}\overline{g}^{\gamma\eta} \partial_\beta \overline{g}_{\xi\eta}(-\frac{1}{2}\partial_\gamma \overline{g}_{ij}+\partial_j\overline{g}_{i\gamma})
\\&=&g^{ij}v^\beta\partial_\beta\overline{\Gamma}_{ij}^\alpha+g^{ij}v^\beta \overline{g}^{\alpha\xi}\Gamma_{ij}^l\partial_\beta \overline{g}_{l\xi}
+g^{ij}v^\beta \overline{g}^{\alpha\xi}\overline{\Gamma}_{ij}^\eta  \partial_\beta \overline{g}_{\xi\eta}.
\end{eqnarray*}
It then follows from (\ref{secondfundamentalform}) and (\ref{vH}) that
\begin{eqnarray*}
\lefteqn{g^{ij}[-\frac{1}{2}\overline{g}^{\alpha\gamma}\partial_\gamma\partial_\beta \overline{g}_{ij}v^\beta
+\overline{g}^{\alpha\gamma}\partial_j\partial_\beta \overline{g}_{i\gamma}v^\beta]}
\\&=&g^{ij}v^\beta\partial_\beta\overline{\Gamma}_{ij}^\alpha+g^{ij}\Gamma_{ij}^l\overline{g}^{\alpha\xi}v^\beta\partial_\beta \overline{g}_{l\xi}
+2n\overline{g}^{\alpha\xi} v^\beta v^\eta\partial_\beta \overline{g}_{\xi\eta}.
\end{eqnarray*}
That is (\ref{Ib}). $\hfill \Box$

\begin{lemma}
\begin{eqnarray}
I^\alpha&=&\triangle^\perp v^\alpha-g^{ij}g^{kl}\langle A_{ik},v\rangle A_{jl}^\alpha+g^{ij}g^{\alpha\gamma}v^\beta\overline{R}_{i\beta j\gamma}
+(4-2n)v^\beta v^\eta \overline{\Gamma}^\alpha_{\beta\eta}. \label{I3}
\end{eqnarray}
\end{lemma}
Proof. It follows from (\ref{I}), (\ref{Ia}) and (\ref{Ib}) that
\begin{eqnarray}\label{I2}
I^\alpha&=&\triangle^\perp v^\alpha-g^{ij}g^{kl}\langle A_{ik},v\rangle A_{jl}^\alpha\nonumber
\\&&-g^{ij}v^\beta[\partial_j\overline{\Gamma}_{i\beta}^\alpha-\partial_\beta\overline{\Gamma}_{ij}^\alpha
+\overline{\Gamma}_{i\beta}^\gamma\overline{\Gamma}_{j\gamma}^\alpha
+\overline{\Gamma}_{i\beta}^k\overline{\Gamma}_{jk}^\alpha-\Gamma_{ij}^k\overline{\Gamma}_{k\beta}^\alpha ]
\\&&+4\overline{g}^{\alpha\gamma}\partial_\beta \overline{g}_{\eta\gamma}v^\beta v^\eta
-2\overline{g}^{\alpha\gamma}\partial_\gamma \overline{g}_{\beta\eta}v^\beta v^\eta.\nonumber
\end{eqnarray}
Note that
\begin{eqnarray}\label{curvature}
\overline{g}^{\alpha\eta}\overline{R}_{j\beta \eta i}
&=&\overline{g}^{\alpha\eta}\langle \overline{\nabla}_j\overline{\nabla}_\beta\partial_i-\overline{\nabla}_\beta\overline{\nabla}_j\partial_i,\partial_\eta\rangle \nonumber
\\&=&\partial_j\overline{\Gamma}_{\beta i}^\alpha+\overline{\Gamma}_{\beta i}^\gamma\overline{\Gamma}_{j\gamma}^\alpha
+\overline{\Gamma}_{\beta i}^k\overline{\Gamma}_{jk}^\alpha
-\partial_\beta \overline{\Gamma}_{ij}^\alpha
-\overline{\Gamma}_{ij}^k \overline{\Gamma}_{\beta k}^\alpha-\overline{\Gamma}_{ij}^\gamma \overline{\Gamma}_{\beta \gamma}^\alpha.
\end{eqnarray}
(\ref{I3}) then follows from (\ref{I2}) and (\ref{curvature}). $\hfill \Box$

\begin{proposition}\label{kgeq3}
Let $Y^{n+1}$, $n\geq 3$, be a minimal surface in $(X^{d+1}, g_+)$ and of the form $Y=(x,u(x,r),r)$ near the boundary of $X$,
where
$$u(x,r)=\frac{1}{2n}H(x)r^2+u^{(4)}(x)r^4+\cdots$$
Then we have
\begin{eqnarray}
8n(n-2)u^{(4)\alpha}&=&\triangle^\perp H^\alpha+g^{ij}g^{kl}\langle A_{ik},H\rangle A_{jl}^\alpha-\frac{2}{n^2}|H|^2H^\alpha\nonumber
\\&&+g^{ij}\overline{g}^{\alpha\gamma}H^\beta\overline{R}_{i\beta j\gamma}+4\overline{g}^{\alpha\gamma}g^{(2)}_{\beta\gamma}H^\beta+
2g^{ij}g^{(2)}_{ij}H^\alpha-2ng^{ik}g^{jl}g^{(2)}_{kl}A_{ij}^\alpha \label{u4}
\\&&-n\overline{g}^{\alpha\gamma}g^{ij}(\overline{\nabla}_\gamma g^{(2)}_{ij}-2\overline{\nabla}_ig^{(2)}_{j\gamma})
-\frac{n-2}{n}H^\beta H^\eta \overline{\Gamma}^\alpha_{\beta\eta}.\nonumber
\end{eqnarray}
\end{proposition}
Proof. It follows from (\ref{w2=}), (\ref{I3}) and (\ref{IIformula}) that
\begin{eqnarray*}
4(n-2)u^{(4)\alpha}&=&-2(n-2)\overline{g}^{\alpha\gamma}g^{(2)}_{\beta\gamma}v^\beta
-8|v|^2v^\alpha+2g^{ij}g^{(2)}_{ij}v^\alpha \nonumber
\\&&-g^{ik}g^{jl}g^{(2)}_{kl}A_{ij}^\alpha+2g^{ik}g^{jl}\langle A_{kl},v\rangle A_{ij}^\alpha \nonumber
\\&&+\triangle^\perp v^\alpha-g^{ij}g^{kl}\langle A_{ik},v\rangle A_{jl}^\alpha+g^{ij}g^{\alpha\gamma}v^\beta\overline{R}_{i\beta j\gamma}
+(4-2n)v^\beta v^\eta \overline{\Gamma}^\alpha_{\beta\eta}
\\&&-\frac{1}{2}\overline{g}^{\alpha\gamma}g^{ij}(\overline{\nabla}_\gamma g^{(2)}_{ij}-2\overline{\nabla}_ig^{(2)}_{j\gamma})
+\overline{g}^{\alpha\gamma}H^\beta g^{(2)}_{\beta\gamma}.
\end{eqnarray*}
Using (\ref{vH}), we get (\ref{u4}). $\hfill \Box$

Substituting $g^{(2)}=-\overline{P}$ and using (\ref{WeylSchouten}), one gets (\ref{u4intro}).
For $n=2$, similar calculations imply the following

\begin{proposition}
Let $Y^3$ be a minimal surface in $(X^{d+1}, g_+)$ and of the form $Y=(x,u(x,r),r)$ near the boundary of $X$,
where
$$u(x,r)=\frac{1}{4}H(x)r^2+w_2(x)r^4\log r +u^{(4)}(x)r^4+\cdots$$
Then we have
\begin{eqnarray}
-16w_2^\alpha(x)&=&\triangle^\perp H^\alpha+g^{ij}g^{kl}\langle A_{ik},H\rangle A_{jl}^\alpha-\frac{1}{2}|H|^2H^\alpha \nonumber
\\&&+g^{ij}\overline{g}^{\alpha\gamma}H^\beta\overline{W}_{i\beta j\gamma}-2\overline{g}^{\alpha\gamma}\overline{P}_{\beta\gamma}H^\beta
-g^{ij}\overline{P}_{ij}H^\alpha+4g^{ik}g^{jl}\overline{P}_{kl}A_{ij}^\alpha \label{w2}
\\&&+2\overline{g}^{\alpha\gamma}g^{ij}(\overline{\nabla}_\gamma \overline{P}_{ij}-2\overline{\nabla}_i\overline{P}_{j\gamma}). \nonumber
\end{eqnarray}
\end{proposition}

\section{Proof of Proposition \ref{GWinvariant4}}

In this section we prove Proposition \ref{GWinvariant4}. Namely, we calculate Graham-Witten's conformal invariant $L_4$ for closed four dimensional
submanifolds. The volume form and the volume expansion of the minimal surface $Y^{n+1}\hookrightarrow (X^{d+1}, g_+)$
are given by \cite{GrahamWitten}
\begin{equation}\label{vkvolumeelement}
d\mu_Y=r^{-n-1}\sqrt{\det \overline{h}}dxdr=r^{-n-1}[v^{(0)}+v^{(2)}r^2+(\mbox{even powers})+v^{(n)}r^n+\cdots]d\mu_\Sigma dr,
\end{equation}
and for $n$ even
\begin{equation}\label{areaexansioineven2}
Vol_{g_+}(Y\cap \{r>\epsilon\})=c_0\epsilon^{-n}+(\mbox{even powers})+c_{n-2}\epsilon^{-2}+L_n\log\frac{1}{\epsilon}+c_n+o(1),
\end{equation}
where
\begin{equation}\label{Ln}
L_n=\int_{\Sigma^n} v^{(n)}d\mu_\Sigma
\end{equation}
is invariant under conformal transformations of $(M^d, \overline{g})$.

Let $n\geq 4$ and
$$u^\alpha(x,r)=v^\alpha(x) r^2+w^\alpha(x) r^4+\cdots$$
where $v(x)=\frac{1}{2n}H$ and  $w=u^{(4)}$ is given by (\ref{u4}). We now compute $v^{(4)}$ in (\ref{vkvolumeelement}). Recall that
$$\overline{h}_{rr}=\overline{h}(Y_r,Y_r)=1+\widetilde{g}_{\alpha\beta}u_r^\alpha u_r^\beta,$$
$$\overline{h}_{ir}=\overline{h}(Y_i,Y_r)=\widetilde{g}_{i\alpha}u_r^\alpha+\widetilde{g}_{\alpha\beta}u_i^\alpha u_r^\beta,$$
and
$$\overline{h}_{ij}=\overline{h}(Y_i,Y_j)=\widetilde{g}_{ij}+\widetilde{g}_{i\alpha}u_j^\alpha
+\widetilde{g}_{j\alpha}u_i^\alpha+\widetilde{g}_{\alpha\beta}u_i^\alpha u_j^\beta.$$
Hence
\begin{eqnarray*}
\overline{h}_{rr}&=&1+u_r^\alpha u_r^\beta \widetilde{g}_{\alpha\beta}
\\&=&1+(2rv^\alpha+4r^3w^\alpha)(2rv^\beta+4r^3w^\beta)(\overline{g}_{\alpha\beta}+g^{(2)}_{\alpha\beta}r^2
+\partial_\gamma \overline{g}_{\alpha\beta}v^\gamma r^2)+o(r^4)
\\&=&1+4r^2|v|^2+4r^4(v^\alpha v^\beta g^{(2)}_{\alpha\beta}+v^\alpha v^\beta v^\gamma \partial_\gamma \overline{g}_{\alpha\beta}+4v^\alpha w^\beta \overline{g}_{\alpha\beta})+o(r^4),
\end{eqnarray*}
\begin{eqnarray*}
\overline{h}_{ir}&=&\widetilde{g}_{i\alpha}u_r^\alpha+u_i^\alpha u_r^\beta \widetilde{g}_{\alpha\beta}
\\&=&2rv^\alpha (g^{(2)}_{i\alpha}r^2+\partial_\gamma \overline{g}_{i\alpha}v^\gamma r^2)+v_i^\alpha r^22v^\beta r\overline{g}_{\alpha\beta}+o(r^4)
\\&=&2r^3v^\alpha( g^{(2)}_{i\alpha}+\partial_\gamma \overline{g}_{i\alpha}v^\gamma +v_i^\beta \overline{g}_{\alpha\beta})+o(r^4),
\end{eqnarray*}
and
\begin{eqnarray*}
\overline{h}_{ij}&=&\widetilde{g}_{ij}+\widetilde{g}_{i\alpha}u_j^\alpha
+\widetilde{g}_{j\alpha}u_i^\alpha+\widetilde{g}_{\alpha\beta}u_i^\alpha u_j^\beta
\\&=&g_{ij}+(g^{(2)}_{ij}+\partial_\gamma \overline{g}_{ij}v^\gamma)r^2
+r^4(g^{(4)}_{ij}+\partial_\beta g^{(2)}_{ij}v^\beta+\partial_\gamma \overline{g}_{ij}w^\gamma+\frac{1}{2}\partial_\beta\partial_\gamma \overline{g}_{ij}v^\beta v^\gamma)
\\&&+r^4(g^{(2)}_{i\alpha}+\partial_\gamma \overline{g}_{i\alpha}v^\gamma)v_j^\alpha+r^4(g^{(2)}_{j\alpha}+\partial_\gamma \overline{g}_{j\alpha}v^\gamma)v_i^\alpha
+r^4\overline{g}_{\alpha\beta}v_i^\alpha v_j^\beta+o(r^4).
\end{eqnarray*}
We rewrite the above as
\begin{eqnarray*}
\overline{h}_{rr}&=&1+ar^2+br^4,
\end{eqnarray*}
\begin{eqnarray*}
\overline{h}_{ir}&=&a_ir^3,
\end{eqnarray*}
\begin{eqnarray*}
\overline{h}_{ij}&=&g_{ij}+a_{ij}r^2+b_{ij}r^4,
\end{eqnarray*}
where
$$a=4|v|^2,$$
$$b=4(v^\alpha v^\beta g^{(2)}_{\alpha\beta}+v^\alpha v^\beta v^\gamma \partial_\gamma \overline{g}_{\alpha\beta}+4\langle v,w\rangle),$$
$$a_i=2v^\alpha( g^{(2)}_{i\alpha}+\partial_\gamma \overline{g}_{i\alpha}v^\gamma +v_i^\beta \overline{g}_{\alpha\beta}),$$
$$a_{ij}=g^{(2)}_{ij}-2\langle A_{ij},v\rangle ,$$
\begin{eqnarray*}
b_{ij}&=&g^{(4)}_{ij}-2\langle A_{ij},w\rangle +\partial_\beta g^{(2)}_{ij}v^\beta+v_i^\alpha g^{(2)}_{j\alpha}+v_j^\alpha g^{(2)}_{i\alpha}
\\&&+\overline{g}_{\alpha\beta}v_i^\alpha v_j^\beta+\frac{1}{2}\partial_\beta\partial_\gamma \overline{g}_{ij}v^\beta v^\gamma
+\partial_\gamma \overline{g}_{i\alpha}v^\gamma v_j^\alpha+\partial_\gamma \overline{g}_{j\alpha}v^\gamma v_i^\alpha.
\end{eqnarray*}
Let $A=(a_{ij}), B=(b_{ij})$. Then
\begin{eqnarray*}
\det \overline{h}&=&(1+ar^2+br^4)\det(\overline{h}_{ij})+o(r^4)
\\&=&\det(\overline{h}_{ij})+ar^2\det(g_{ij}+a_{ij}r^2)+br^4\det(g_{ij})+o(r^4)
\\&=&\det(\overline{h}_{ij})+ar^2\det(g_{ij})(1+g^{ij}a_{ij}r^2)+br^4\det(g_{ij})+o(r^4).
\end{eqnarray*}
Note that
\begin{eqnarray*}
\det(\overline{h}_{ij})&=&\det(g_{ij})[1+g^{ij}a_{ij}r^2+(g^{ij}b_{ij}+\sigma_2(g^{-1}A))r^4]+o(r^4).
\end{eqnarray*}
Therefore,
\begin{eqnarray*}
\det \overline{h}&=&\det(g_{ij})(1+Cr^2+Dr^4)+o(r^4),
\end{eqnarray*}
where
$$C=a+g^{ij}a_{ij},$$
$$D=atr(g^{-1}A)+\sigma_2(g^{-1}A)+b+tr(g^{-1}B).$$
Then
\begin{eqnarray}
d\mu_Y&=&r^{-(n+1)}\sqrt{\det \overline{h}}dxdr \nonumber
\\&=&r^{-(n+1)} (1+\frac{C}{2}r^2+\frac{1}{2}(D-\frac{C^2}{4})r^4+o(r^4))d\mu_\Sigma dr, \label{volumeformofY}
\end{eqnarray}
and the coefficient $v^{(4)}$ in (\ref{vkvolumeelement}) is
\begin{equation}\label{v4}
v^{(4)}=\frac{1}{2}(D-\frac{C^2}{4}).
\end{equation}

\begin{proposition}
Let $Y^{n+1}$, $n\geq 4$, be a minimal surface in $(X^{d+1}, g_+)$.
Then we have
\begin{eqnarray}
2v^{(4)}&=&4|v|^2(g^{ij}g^{(2)}_{ij}-4n|v|^2)+\frac{1}{2}(g^{ij}g^{(2)}_{ij}-4n|v|^2)^2-\frac{1}{2}|g^{(2)}_{ij}-2\langle A_{ij},v\rangle |_g^2 \nonumber
\\&&-\frac{1}{4}(g^{ij}g^{(2)}_{ij}-4(n-1)|v|^2)^2+4g^{(2)}_{\alpha\beta}v^\alpha v^\beta -4(n-4)\langle v,w\rangle+g^{ij}g^{(4)}_{ij}  \label{v40}
\\&&+I+II, \label{v4I-II}\nonumber
\end{eqnarray}
where
\begin{equation}\label{I41}
I:=g^{ij}\overline{g}_{\alpha\beta}v_i^\alpha v_j^\beta
+\frac{1}{2}g^{ij}\partial_\beta\partial_\gamma \overline{g}_{ij}v^\beta v^\gamma+2g^{ij}\partial_\gamma \overline{g}_{i\alpha}v^\gamma v_j^\alpha+4v^\alpha v^\beta v^\gamma\partial_\gamma \overline{g}_{\alpha\beta},
\end{equation}
and
\begin{eqnarray}\label{II41}
II&:=&g^{ij}(\partial_\beta g^{(2)}_{ij}v^\beta+2g^{(2)}_{i\alpha}v_j^\alpha).
\end{eqnarray}
\end{proposition}

We first deal with (\ref{II41}).

\begin{lemma} We have
\begin{eqnarray}
II&=&g^{ij}v^\beta(\overline{\nabla}_\beta g^{(2)}_{ij}-2\overline{\nabla}_jg^{(2)}_{i\beta})-4ng^{(2)}_{\alpha\beta}v^\alpha v^\beta
+2g^{ij}\nabla_j(v^\beta g^{(2)}_{i\beta}). \label{II42}
\end{eqnarray}
\end{lemma}
Proof. Note that
\begin{eqnarray*}
g^{ij}\nabla_j(v^\beta g^{(2)}_{i\beta})&=&g^{ij}[\partial_j(v^\beta g^{(2)}_{i\beta})-\Gamma_{ij}^kv^\beta g^{(2)}_{k\beta}]
\\&=&g^{ij}v_j^\beta g^{(2)}_{i\beta}+g^{ij}v^\beta \partial_jg^{(2)}_{i\beta}-g^{ij}\Gamma_{ij}^kv^\beta g^{(2)}_{k\beta},
\end{eqnarray*}
hence
\begin{eqnarray*}
\lefteqn{g^{ij}(\partial_\beta g^{(2)}_{ij}v^\beta+2g^{(2)}_{i\beta}v_j^\beta)-2g^{ij}\nabla_j(v^\beta g^{(2)}_{i\beta})}
\\&=&g^{ij}v^\beta(\partial_\beta g^{(2)}_{ij}-2\partial_jg^{(2)}_{i\beta}+2\Gamma_{ij}^kg^{(2)}_{k\beta})
\\&=&g^{ij}v^\beta(\partial_\beta g^{(2)}_{ij}-\partial_jg^{(2)}_{i\beta}-\partial_ig^{(2)}_{j\beta}+2\Gamma_{ij}^kg^{(2)}_{k\beta}).
\end{eqnarray*}
As we have shown in (\ref{II}), we have
\begin{eqnarray*}
\partial_\beta g^{(2)}_{ij}-\partial_jg^{(2)}_{i\beta}-\partial_ig^{(2)}_{j\beta}
&=&\overline{\nabla}_\beta g^{(2)}_{ij}-\overline{\nabla}_jg^{(2)}_{i\beta}-\overline{\nabla}_ig^{(2)}_{j\beta}-2\overline{\Gamma}_{ij}^Dg^{(2)}_{D\beta}.
\end{eqnarray*}
Therefore,
\begin{eqnarray*}
\lefteqn{g^{ij}(\partial_\beta g^{(2)}_{ij}v^\beta+2g^{(2)}_{i\beta}v_j^\beta)-2g^{ij}\nabla_j(v^\beta g^{(2)}_{i\beta})}
\\&=&g^{ij}v^\beta(\overline{\nabla}_\beta g^{(2)}_{ij}-2\overline{\nabla}_jg^{(2)}_{i\beta}
-2\overline{\Gamma}_{ij}^Dg^{(2)}_{D\beta}+2\Gamma_{ij}^kg^{(2)}_{k\beta})
\\&=&g^{ij}v^\beta(\overline{\nabla}_\beta g^{(2)}_{ij}-2\overline{\nabla}_jg^{(2)}_{i\beta}-2\overline{\Gamma}_{ij}^\alpha g^{(2)}_{\alpha\beta})
\\&=&g^{ij}v^\beta(\overline{\nabla}_\beta g^{(2)}_{ij}-2\overline{\nabla}_jg^{(2)}_{i\beta})-2H^\alpha v^\beta g^{(2)}_{\alpha\beta}.
\end{eqnarray*}
$\hfill \Box$

\begin{lemma}
We have
\begin{eqnarray}
I&=&-\langle \triangle^\perp v,v\rangle +g^{ij}g^{kl}\langle A_{ik},v\rangle \langle A_{jl},v\rangle
-g^{ij}v^\alpha v^\beta\overline{R}_{i\alpha j\beta}  \label{I43}
\\&&-2(n-4)\overline{g}_{\alpha\xi}v^\xi v^\beta v^\gamma \overline{\Gamma}_{\beta\gamma}^\alpha
+g^{ij}\nabla_j(\overline{g}_{\alpha\beta}v_i^\alpha v^\beta+\partial_\alpha \overline{g}_{i\beta} v^\alpha v^\beta).\nonumber
\end{eqnarray}
\end{lemma}
Proof. Note that
\begin{eqnarray*}
g^{ij}\overline{g}_{\alpha\beta}v_i^\alpha v_j^\beta&=&g^{ij}\nabla_j(\overline{g}_{\alpha\beta}v_i^\alpha v^\beta)
+g^{ij}\Gamma_{ij}^k\overline{g}_{\alpha\beta}v_k^\alpha v^\beta
-g^{ij}\partial_j\overline{g}_{\alpha\beta}v_i^\alpha v^\beta-g^{ij}\overline{g}_{\alpha\beta}\partial_jv_i^\alpha v^\beta.
\end{eqnarray*}
It then follows from (\ref{Ia}):
\begin{eqnarray*}
\partial_jv_i^\alpha&=&\nabla_j^\perp\nabla_i^\perp v^\alpha-g^{kl}\langle A_{ik},v\rangle A_{jl}^\alpha
-v_i^\beta \overline{\Gamma}_{j\beta}^\alpha-v_j^\beta \overline{\Gamma}_{i\beta}^\alpha+v_k^\alpha\Gamma_{ij}^k \nonumber
\\&&-v^\beta(\partial_j\overline{\Gamma}_{i\beta}^\alpha+\overline{\Gamma}_{i\beta}^\gamma\overline{\Gamma}_{j\gamma}^\alpha
+\overline{\Gamma}_{i\beta}^k\overline{\Gamma}_{jk}^\alpha-\Gamma_{ij}^k\overline{\Gamma}_{k\beta}^\alpha)
\end{eqnarray*}
and (\ref{I41}):
$$I=g^{ij}\overline{g}_{\alpha\beta}v_i^\alpha v_j^\beta
+\frac{1}{2}g^{ij}\partial_\beta\partial_\gamma \overline{g}_{ij}v^\beta v^\gamma+2g^{ij}\partial_\gamma \overline{g}_{i\alpha}v^\gamma v_j^\alpha+4v^\alpha v^\beta v^\gamma\partial_\gamma \overline{g}_{\alpha\beta},$$
that
\begin{eqnarray*}
I&=&-\langle\triangle^\perp v,v\rangle+g^{ij}g^{kl}\langle A_{ik},v\rangle \langle A_{jl},v\rangle
+g^{ij}\nabla_j(\overline{g}_{\alpha\beta}v_i^\alpha v^\beta)
\\&&+\frac{1}{2}g^{ij}\partial_\beta\partial_\gamma \overline{g}_{ij}v^\beta v^\gamma+2g^{ij}\partial_\gamma \overline{g}_{i\alpha}v^\gamma v_j^\alpha
+4v^\alpha v^\beta v^\gamma\partial_\gamma \overline{g}_{\alpha\beta}-g^{ij}\partial_j\overline{g}_{\alpha\beta}v_i^\alpha v^\beta
\\&&+g^{ij}\overline{g}_{\alpha\xi}v^\xi [2v_i^\beta \overline{\Gamma}_{j\beta}^\alpha
+v^\beta(\partial_j\overline{\Gamma}_{i\beta}^\alpha+\overline{\Gamma}_{i\beta}^\gamma\overline{\Gamma}_{j\gamma}^\alpha
+\overline{\Gamma}_{i\beta}^k\overline{\Gamma}_{jk}^\alpha-\Gamma_{ij}^k\overline{\Gamma}_{k\beta}^\alpha)]
\end{eqnarray*}
Note that
\begin{eqnarray*}
\lefteqn{2g^{ij}\partial_\gamma \overline{g}_{i\alpha}v^\gamma v_j^\alpha-g^{ij}\partial_j\overline{g}_{\alpha\beta}v_i^\alpha v^\beta
+2g^{ij}\overline{g}_{\alpha\xi}v^\xi v_i^\beta \overline{\Gamma}_{j\beta}^\alpha}
\\&=&g^{ij}v_i^\alpha v^\gamma[2\partial_\gamma \overline{g}_{j\alpha}-\partial_j\overline{g}_{\alpha\gamma}+2\overline{g}_{\beta\gamma}\overline{\Gamma}_{j\alpha}^\beta]
\\&=&g^{ij}v_i^\alpha v^\gamma[\partial_\gamma \overline{g}_{j\alpha}+\partial_\alpha \overline{g}_{j\gamma}]
\\&=&g^{ij}\partial_\gamma \overline{g}_{j\alpha}\partial_i(v^\alpha v^\gamma)
\\&=&g^{ij}\nabla_i(\partial_\gamma \overline{g}_{j\alpha} v^\alpha v^\gamma)
+g^{ij}\Gamma_{ij}^k\partial_\gamma \overline{g}_{k\alpha} v^\alpha v^\gamma-g^{ij}\partial_i\partial_\gamma \overline{g}_{j\alpha}v^\alpha v^\gamma.
\end{eqnarray*}
Hence
\begin{eqnarray*}
I&=&-\langle\triangle^\perp v,v\rangle+g^{ij}g^{kl}\langle A_{ik},v\rangle \langle A_{jl},v\rangle
+g^{ij}\nabla_j(\overline{g}_{\alpha\beta}v_i^\alpha v^\beta+\partial_\gamma \overline{g}_{i\alpha} v^\alpha v^\gamma)
\\&&+\frac{1}{2}g^{ij}\partial_\beta\partial_\gamma \overline{g}_{ij}v^\beta v^\gamma-g^{ij}\partial_j\partial_\gamma \overline{g}_{i\beta}v^\beta v^\gamma
+g^{ij}\Gamma_{ij}^k\partial_\gamma \overline{g}_{k\alpha} v^\alpha v^\gamma+4v^\alpha v^\beta v^\gamma\partial_\gamma \overline{g}_{\alpha\beta}
\\&&+g^{ij}\overline{g}_{\alpha\xi}v^\xi v^\beta(\partial_j\overline{\Gamma}_{i\beta}^\alpha+\overline{\Gamma}_{i\beta}^\gamma\overline{\Gamma}_{j\gamma}^\alpha
+\overline{\Gamma}_{i\beta}^k\overline{\Gamma}_{jk}^\alpha-\Gamma_{ij}^k\overline{\Gamma}_{k\beta}^\alpha).
\end{eqnarray*}

Recall (\ref{Ib}):
\begin{eqnarray*}
\lefteqn{g^{ij}[-\frac{1}{2}\overline{g}^{\alpha\gamma}\partial_\gamma\partial_\beta \overline{g}_{ij}v^\beta
+\overline{g}^{\alpha\gamma}\partial_j\partial_\beta \overline{g}_{i\gamma}v^\beta]}\nonumber
\\&=&g^{ij}v^\beta\partial_\beta\overline{\Gamma}_{ij}^\alpha+g^{ij}\Gamma_{ij}^l\overline{g}^{\alpha\xi}v^\beta\partial_\beta \overline{g}_{l\xi}
+2n\overline{g}^{\alpha\xi} v^\beta v^\eta  \partial_\beta \overline{g}_{\xi\eta},
\end{eqnarray*}
hence
\begin{eqnarray*}
\lefteqn{\frac{1}{2}g^{ij}\partial_\beta\partial_\gamma \overline{g}_{ij}v^\beta v^\gamma
-g^{ij}\partial_j\partial_\gamma \overline{g}_{i\beta}v^\beta v^\gamma}
\\&=&- \overline{g}_{\alpha\zeta}v^\zeta g^{ij}[-\frac{1}{2}\overline{g}^{\alpha\gamma}\partial_\gamma\partial_\beta \overline{g}_{ij}v^\beta
+\overline{g}^{\alpha\gamma}\partial_j\partial_\beta \overline{g}_{i\gamma}v^\beta]
\\&=&-g^{ij}\overline{g}_{\alpha\xi}v^\xi v^\beta \partial_\beta\overline{\Gamma}_{ij}^\alpha-g^{ij}\Gamma_{ij}^lv^\xi v^\beta \partial_\beta \overline{g}_{l\xi}
-2n v^\xi v^\beta v^\eta\partial_\beta \overline{g}_{\xi\eta}.
\end{eqnarray*}
Then
\begin{eqnarray*}
I&=&-\langle \triangle^\perp v,v\rangle+g^{ij}g^{kl}\langle A_{ik},v\rangle \langle A_{jl},v\rangle
+g^{ij}\nabla_j(\overline{g}_{\alpha\beta}v_i^\alpha v^\beta+\partial_\gamma \overline{g}_{i\alpha} v^\alpha v^\gamma)
\\&&+g^{ij}\overline{g}_{\alpha\xi}v^\xi v^\beta(\partial_j\overline{\Gamma}_{i\beta}^\alpha
+\overline{\Gamma}_{i\beta}^\gamma\overline{\Gamma}_{j\gamma}^\alpha
+\overline{\Gamma}_{i\beta}^k\overline{\Gamma}_{jk}^\alpha-\partial_\beta\overline{\Gamma}_{ij}^\alpha-\Gamma_{ij}^k\overline{\Gamma}_{k\beta}^\alpha)
\\&&+(4-2n)v^\alpha v^\beta v^\gamma\partial_\gamma \overline{g}_{\alpha\beta}.
\end{eqnarray*}
Recall (\ref{curvature}):
\begin{eqnarray*}
\overline{g}^{\alpha\eta}\overline{R}_{j\beta \eta i}
&=&\partial_j\overline{\Gamma}_{\beta i}^\alpha+\overline{\Gamma}_{\beta i}^\gamma\overline{\Gamma}_{j\gamma}^\alpha
+\overline{\Gamma}_{\beta i}^k\overline{\Gamma}_{jk}^\alpha
-\partial_\beta \overline{\Gamma}_{ij}^\alpha
-\overline{\Gamma}_{ij}^k \overline{\Gamma}_{\beta k}^\alpha-\overline{\Gamma}_{ij}^\gamma \overline{\Gamma}_{\beta \gamma}^\alpha,
\end{eqnarray*}
and $\overline{\Gamma}_{ij}^\gamma=A_{ij}^\gamma$,  $H^\gamma=2nv^\gamma$, so we have
\begin{eqnarray*}
I&=&-\langle \triangle^\perp v,v\rangle +g^{ij}g^{kl}\langle A_{ik},v\rangle \langle A_{jl},v\rangle
+g^{ij}\nabla_j(\overline{g}_{\alpha\beta}v_i^\alpha v^\beta+\partial_\gamma \overline{g}_{i\alpha} v^\alpha v^\gamma)
\\&&-g^{ij}v^\alpha v^\beta\overline{R}_{i\alpha j\beta}-2(n-4)\overline{g}_{\alpha\xi}v^\xi v^\beta v^\gamma \overline{\Gamma}_{\beta\gamma}^\alpha.
\end{eqnarray*}
$\hfill \Box$

It follows from (\ref{v40}), (\ref{I43}) and (\ref{II42}) the following

\begin{proposition}
Let $Y^{n+1}$, $n\geq 4$, be a minimal surface in $(X^{d+1}, g_+)$. The coefficient $v^{(4)}$ in (\ref{vkvolumeelement}) is given by
\begin{eqnarray}
2v^{(4)}&=&-\langle\triangle^\perp v,v\rangle-g^{ik}g^{jl}\langle A_{ij},v\rangle \langle A_{kl},v\rangle
+4(n^2-2n-1)|v|^4\nonumber
\\&&+\frac{1}{4}(g^{ij}g^{(2)}_{ij})^2-\frac{1}{2}|g^{(2)}_{ij}|_g^2 +g^{ij}g^{(4)}_{ij}\nonumber
\\&&+2g^{ik}g^{jl}g^{(2)}_{ij}\langle A_{kl},v\rangle-2(n-1)g^{ij}g^{(2)}_{ij}|v|^2
-4(n-1) g^{(2)}_{\alpha\beta}v^\alpha v^\beta  -g^{ij}\overline{R}_{i\alpha j\beta}v^\alpha v^\beta \nonumber
\\&&+g^{ij}v^\beta(\overline{\nabla}_\beta g^{(2)}_{ij}-2\overline{\nabla}_ig^{(2)}_{j\beta}) \label{v4prop}
\\&&-4(n-4)\langle v,w\rangle-2(n-4)\overline{g}_{\alpha\xi}v^\xi v^\beta v^\gamma \overline{\Gamma}_{\beta\gamma}^\alpha  \nonumber
\\&&+g^{ij}\nabla_j(\overline{g}_{\alpha\beta}v_i^\alpha v^\beta+\partial_\alpha \overline{g}_{i\beta} v^\alpha v^\beta+2g^{(2)}_{i\alpha}v^\alpha), \nonumber
\end{eqnarray}
where $v=\frac{1}{2n}H$ and $w=u^{(4)}$ is given by (\ref{u4}).
\end{proposition}

Notice that the second last line in (\ref{v4prop}) vanishes for $n=4$ and for $n\geq 5$ the part involving $\overline{\Gamma}$ cancels.
Note also that the last line is a divergence which is independent of the choice of local coordinates $(x^i, y^\alpha)$.
By integrating (\ref{v4prop}) over a closed submanifold $\Sigma^4$ and using (\ref{g2}), (\ref{WeylSchouten}) and (\ref{g4}),
we obtain the expression (\ref{L4full}) of Graham-Witten's conformal invariant $L_4$.

\section{$L_4$ for closed submanifolds of Euclidean spaces}

In this section, we assume $n=4, d\geq 5$ and the ambient space is $(M^d, \overline{g})=\mathbb{R}^d$. It follows from (\ref{L4full}) that
\begin{eqnarray*}
128L_4&=&\int_\Sigma (|\nabla^\perp H|^2-g^{ik}g^{jl}\langle A_{ij},H\rangle \langle A_{kl},H\rangle +\frac{7}{16}|H|^4)d\mu_g.
\end{eqnarray*}
We now consider the following functional for closed four dimensional submanifolds $F: \Sigma\rightarrow \mathbb{R}^d$
\begin{equation*}
\mathcal{L}_4(\Sigma^4)=\frac{1}{2}\int_\Sigma (|\nabla^\perp H|^2-g^{ik}g^{jl}\langle A_{ij},H\rangle \langle A_{kl},H\rangle +\frac{7}{16}|H|^4)d\mu_g.
\end{equation*}
In the sequel we use local coordinates $(x^i)$ of $\Sigma$  which is normal at the point of consideration with respect to the induced metric
$$g_{ij}:=\langle F_i, F_j\rangle.$$

\begin{proposition}
$F:\Sigma^4\hookrightarrow \mathbb{R}^d$ is a critical point of $\mathcal{L}_4$ if and only if
\begin{equation}\label{EulerLagrange}
\mathcal{E}=0,
\end{equation}
where
\begin{eqnarray}
\mathcal{E}&=&\triangle^\perp(\triangle^\perp H+\langle A_{ij},H\rangle A_{ij}-\frac{7}{8}|H|^2H)\nonumber
\\&&+\langle \triangle^\perp H,A_{ij}\rangle A_{ij}+<\triangle^\perp H, H>H-\langle \nabla_i^\perp H,\nabla_j^\perp H\rangle A_{ij}
+\frac{3}{2}|\nabla^\perp H|^2H \nonumber
\\&&+2\langle A_{ij},H\rangle \triangle^\perp A_{ij}+2\langle A_{ij},\nabla_i^\perp H\rangle \nabla_j^\perp H
+3\langle H,\nabla_i^\perp H\rangle \nabla_i^\perp H \nonumber
\\&&-\langle A_{ij},H\rangle \langle A_{jk},H\rangle A_{ik}
-\frac{1}{2}\langle A_{ij},H\rangle \langle A_{ij},H\rangle H+3\langle A_{ij},H\rangle \langle A_{ij},A_{kl}\rangle A_{kl} \label{E}
\\&&-\frac{7}{8}|H|^2\langle A_{ij},H\rangle A_{ij}+\frac{7}{32}|H|^4H \nonumber
\\&&-2\langle A_{ij},H\rangle\langle A_{ik},A_{jl}\rangle A_{kl}+2\langle A_{ij},H\rangle\langle A_{ik}, A_{kl}\rangle A_{jl}.\nonumber
\end{eqnarray}
\end{proposition}
Proof. Let $X\in T^\perp \Sigma$ be a variation field. We have
$$\delta g_{ij}=-2\langle A_{ij},X\rangle,$$
\begin{eqnarray*}
\delta A_{ij}=\nabla_i^\perp \nabla_j^\perp X-\langle A_{jk}, X\rangle A_{ik}-\langle A_{ij},\nabla_k^\perp X\rangle F_k,
\end{eqnarray*}
$$\delta H=\triangle^\perp X+\langle A_{kl},X\rangle A_{kl}-\langle H,\nabla_k^\perp X\rangle F_k,$$
and
$$[\delta \nabla_i^\perp H]^\perp=\nabla_i^\perp(\triangle^\perp X+\langle A_{kl},X\rangle A_{kl})-\langle H,\nabla_k^\perp X\rangle A_{ik}
+\langle A_{ik},H\rangle \nabla_k^\perp X.$$
Hence
\begin{eqnarray*}
\lefteqn{\delta\int_\Sigma |\nabla^\perp H|_g^2d\mu_g}
\\&=&\int_\Sigma(2\langle A_{ij},X\rangle \langle  \nabla_i^\perp H, \nabla_j^\perp H\rangle
-|\nabla^\perp H|^2\langle H,X\rangle )d\mu
\\&&+\int_\Sigma 2\langle \nabla_i^\perp(\triangle^\perp X+\langle A_{kl},X\rangle A_{kl})-\langle H,\nabla_k^\perp X\rangle A_{ik}
+\langle A_{ik},H\rangle \nabla_k^\perp X,\nabla_i^\perp H\rangle d\mu
\\&=&\int_\Sigma \langle -2(\triangle^\perp)^2H-2\langle \triangle^\perp H,A_{ij}\rangle A_{ij}
+2\langle \nabla_i^\perp H, \nabla_j^\perp H\rangle A_{ij}-|\nabla^\perp H|^2H, X\rangle d\mu
\\&&+\int_\Sigma 2\langle  \nabla_k^\perp (\langle A_{ik},\nabla_i^\perp H\rangle H-\langle A_{ik},H\rangle \nabla_i^\perp H),X\rangle d\mu.
\end{eqnarray*}
We then compute
\begin{eqnarray*}
\lefteqn{\delta\int_\Sigma(-g^{ik}g^{jl}\langle A_{ij},H\rangle \langle A_{kl},H\rangle)d\mu_g}
\\&=&\int_\Sigma \langle -4\langle A_{ij}, H\rangle \langle A_{jk}, H\rangle A_{ik}
+\langle A_{ij}, H\rangle \langle A_{ij}, H\rangle H, X\rangle d\mu
\\&&-2\int_\Sigma \langle A_{ij},H\rangle \langle\nabla_i^\perp\nabla_j^\perp X-\langle A_{jk},X\rangle A_{ik},H\rangle d\mu
\\&&-2\int_\Sigma \langle A_{ij},H\rangle \langle A_{ij}, \triangle^\perp X+\langle A_{kl}, X\rangle A_{kl}\rangle d\mu
\\&=&\int_\Sigma \langle -2\nabla_j^\perp \nabla_i^\perp (\langle A_{ij},H\rangle H)-2\triangle^\perp (\langle A_{ij}, H\rangle A_{ij}), X\rangle d\mu
\\&&+\int_\Sigma \langle-2\langle A_{ij}, H\rangle \langle A_{jk}, H\rangle A_{ik}
+\langle A_{ij}, H\rangle \langle A_{ij}, H\rangle H-2\langle A_{ij}, H\rangle \langle A_{ij}, A_{kl}\rangle A_{kl}, X\rangle d\mu,
\end{eqnarray*}
and
\begin{eqnarray*}
\delta\int_\Sigma\frac{7}{16}|H|^4d\mu_g&=&\int_\Sigma\frac{7}{4}|H|^2\langle H,\triangle^\perp X+\langle A_{ij},X\rangle A_{ij}\rangle d\mu
-\int_\Sigma \frac{7}{16}|H|^4\langle H,X\rangle d\mu
\\&=&\int_\Sigma\langle \frac{7}{4}\triangle^\perp(|H|^2H)+\frac{7}{4}|H|^2\langle A_{ij},H\rangle A_{ij}-\frac{7}{16}|H|^4H, X\rangle d\mu.
\end{eqnarray*}
Therefore
\begin{eqnarray}\label{variation}
\delta \mathcal{L}_4(X)=-\int_\Sigma \langle\mathcal{E}, X\rangle d\mu,
\end{eqnarray}
where
\begin{eqnarray}
\mathcal{E}&=&(\triangle^\perp)^2H+\langle \triangle^\perp H,A_{ij}\rangle A_{ij}-\langle \nabla_i^\perp H,\nabla_j^\perp H\rangle A_{ij}+\frac{1}{2}|\nabla^\perp H|^2H \nonumber
\\&&-\nabla_k^\perp (\langle A_{ik},\nabla_i^\perp H\rangle H)+\nabla_k^\perp (\langle A_{ik}, H\rangle \nabla_i^\perp H) \label{GR}
\\&&+\nabla_i^\perp\nabla_j^\perp(\langle A_{ij},H\rangle H)+\triangle^\perp(\langle A_{ij},H\rangle A_{ij})-\frac{7}{8}\triangle^\perp(|H|^2H) \nonumber
\\&&+\langle A_{ij},H\rangle \langle A_{jk},H\rangle A_{ik}
-\frac{1}{2}\langle A_{ij},H\rangle \langle A_{ij},H\rangle H+\langle A_{ij},H\rangle \langle A_{ij},A_{kl}\rangle A_{kl} \nonumber
\\&&-\frac{7}{8}|H|^2\langle A_{ij},H\rangle A_{ij}+\frac{7}{32}|H|^4H. \nonumber
\end{eqnarray}
Note that
\begin{eqnarray*}
\lefteqn{-\nabla_k^\perp (\langle A_{ik},\nabla_i^\perp H\rangle H)+\nabla_k^\perp (\langle A_{ik}, H\rangle \nabla_i^\perp H)
+\nabla_i^\perp\nabla_j^\perp(\langle A_{ij},H\rangle H)}
\\&=&2\langle A_{ij},H\rangle \nabla_i^\perp\nabla_j^\perp H+2\langle A_{ij},\nabla_i^\perp H\rangle \nabla_j^\perp H
+3\langle H,\nabla_i^\perp H\rangle \nabla_i^\perp H
\\&&+\langle \triangle^\perp H, H\rangle H+|\nabla^\perp H|^2H,
\end{eqnarray*}
and
\begin{eqnarray*}
\nabla_i^\perp\nabla_j^\perp H&=&\triangle^\perp A_{ij}+\langle A_{ij}, A_{kl}\rangle A_{kl}-\langle A_{ik},A_{jl}\rangle A_{kl}
+\langle A_{ik}, A_{kl}\rangle A_{jl}-\langle A_{ik}, H\rangle A_{jk},
\end{eqnarray*}
hence
\begin{eqnarray*}
\mathcal{E}&=&\triangle^\perp(\triangle^\perp H+\langle A_{ij},H\rangle A_{ij}-\frac{7}{8}|H|^2H)
\\&&+\langle \triangle^\perp H,A_{ij}\rangle A_{ij}+<\triangle^\perp H, H>H-\langle \nabla_i^\perp H,\nabla_j^\perp H\rangle A_{ij}
+\frac{3}{2}|\nabla^\perp H|^2H
\\&&+2\langle A_{ij},H\rangle \triangle^\perp A_{ij}+2\langle A_{ij},\nabla_i^\perp H\rangle \nabla_j^\perp H
+3\langle H,\nabla_i^\perp H\rangle \nabla_i^\perp H
\\&&-\langle A_{ij},H\rangle \langle A_{jk},H\rangle A_{ik}
-\frac{1}{2}\langle A_{ij},H\rangle \langle A_{ij},H\rangle H+3\langle A_{ij},H\rangle \langle A_{ij},A_{kl}\rangle A_{kl}
\\&&-\frac{7}{8}|H|^2\langle A_{ij},H\rangle A_{ij}+\frac{7}{32}|H|^4H
\\&&+2\langle A_{ij},H\rangle(-\langle A_{ik},A_{jl}\rangle A_{kl}+\langle A_{ik}, A_{kl}\rangle A_{jl}).
\end{eqnarray*}
$\hfill \Box$

We search for critical points of $\mathcal{L}_4$ of the form
\begin{equation}\label{productcriticalpoint}
\mathbb{S}^{k_1}(r_1)\times \cdots \times \mathbb{S}^{k_m}(r_m),
\end{equation}
where $r_i>0, i=1,\cdots, m$, is the radius of a round sphere $\mathbb{S}^{k_i}\hookrightarrow \mathbb{R}^{k_i+1}$, and $k_1\geq k_2\geq \cdots \geq k_m$,
$k_1+\cdots+k_m=4$. Submanifolds of the form (\ref{productcriticalpoint}) have parallel second fundamental form and parallel mean curvature in the normal bundle. We fix $r_1=1$.

For $m=1$, using the equation (\ref{E}) we have the critical point of the round sphere and
$$\mathcal{L}_4(\mathbb{S}^4)=24Vol(\mathbb{S}^4)=64\pi^2.$$

For $m=2$ and $k_1=3$, that is to consider $\mathbb{S}^3(1)\times \mathbb{S}^1(r_2)$, in the case we have
$$A=\delta_{11}\nu_1+\delta_{22}\nu_1+\delta_{33}\nu_1+\frac{1}{r_2}\delta_{44}\nu_2, \quad
H=3\nu_1+\frac{1}{r_2}\nu_2,$$
where $\nu_1$ and $\nu_2$ are unit inner normal vector fields. The equation (\ref{E}) now reads
$$\mathcal{E}=\frac{27}{32}(\frac{1}{r_2^2}-3)(\frac{1}{r_2^2}-\frac{5}{3})(-\nu_1+\frac{1}{r_2}\nu_2),$$
hence we have the following two solutions of the form (\ref{productcriticalpoint}) to (\ref{EulerLagrange})
$$\mathbb{S}^3(1)\times \mathbb{S}^1(\frac{1}{\sqrt{3}}), \quad \mathbb{S}^3(1)\times \mathbb{S}^1(\sqrt{\frac{3}{5}}).$$
The value
$$\mathcal{L}_4[\mathbb{S}^3(1)\times \mathbb{S}^1(r_2)]=2\pi^3r_2[\frac{7}{16}(9+\frac{1}{r_2^2})^2
-27-\frac{1}{r_2^4}].$$
In particular
$$\lim_{r_2\rightarrow +\infty}\mathcal{L}_4[\mathbb{S}^3(1)\times \mathbb{S}^1(r_2)]=+\infty, \quad
\lim_{r_2\rightarrow 0}\mathcal{L}_4[\mathbb{S}^3(1)\times \mathbb{S}^1(r_2)]=-\infty,$$
and
$$\mathcal{L}_4[\mathbb{S}^3(1)\times \mathbb{S}^1(\frac{1}{\sqrt{3}})]=18\sqrt{3}\pi^3,
\quad \mathcal{L}_4[\mathbb{S}^3(1)\times \mathbb{S}^1(\sqrt{\frac{3}{5}})]=8\sqrt{15}\pi^3.$$

For $m=2$ and $k_1=2$, that is to consider $\mathbb{S}^2(1)\times \mathbb{S}^2(r_2)$, in the case we have
$$A=\delta_{11}\nu_1+\delta_{22}\nu_1+\frac{1}{r_2}\delta_{33}\nu_2+\frac{1}{r_2}\delta_{44}\nu_2,
\quad H=2\nu_1+\frac{2}{r_2}\nu_2.$$
The equation (\ref{E}) now reads
$$\mathcal{E}=(\frac{1}{r_2^4}-1)(-\nu_1+\frac{1}{r_2}\nu_2),$$
hence we have the following solution of the form (\ref{productcriticalpoint}) to (\ref{EulerLagrange})
$$\mathbb{S}^2(1)\times \mathbb{S}^2(1).$$
The value
$$\mathcal{L}_4[\mathbb{S}^2(1)\times \mathbb{S}^2(r_2)]=8\pi^2r_2[7(1+\frac{1}{r_2^2})^2
-8(1+\frac{1}{r_2^4})].$$
In particular
$$\lim_{r_2\rightarrow +\infty}\mathcal{L}_4[\mathbb{S}^2(1)\times \mathbb{S}^2(r_2)]=
\lim_{r_2\rightarrow 0}\mathcal{L}_4[\mathbb{S}^2(1)\times \mathbb{S}^2(r_2)]=-\infty$$
and
$$\mathcal{L}_4[\mathbb{S}^2(1)\times \mathbb{S}^2(1)]=96\pi^2.$$

For $m=3$, that is to consider $\mathbb{S}^2(1)\times \mathbb{S}^1(r_2)\times \mathbb{S}^1(r_3)$, in the case we have
$$A=\delta_{11}\nu_1+\delta_{22}\nu_1+\frac{1}{r_2}\delta_{33}\nu_2+\frac{1}{r_3}\delta_{44}\nu_3, \quad
H=2\nu_1+\frac{1}{r_2}\nu_2+\frac{1}{r_3}\nu_3.$$
Let $x=\frac{1}{r_2^2}$ and $y=\frac{1}{r_3^2}$. The equation (\ref{E}) now reads
\begin{eqnarray*}
\mathcal{E}&=&[16-(8+x^2+y^2)-\frac{7}{2}(4+x+y)+\frac{7}{16}(4+x+y)^2]\nu_1
\\&&+[2x^2-\frac{1}{2}(8+x^2+y^2)-\frac{7}{8}(4+x+y)x+\frac{7}{32}(4+x+y)^2]\frac{\nu_2}{r_2}
\\&&+[2y^2-\frac{1}{2}(8+x^2+y^2)-\frac{7}{8}(4+x+y)y+\frac{7}{32}(4+x+y)^2]\frac{\nu_3}{r_3}.
\end{eqnarray*}
Then $\mathcal{E}=0$ if $x=y=2$, or $x=2, y=\frac{10}{9}$.
We have the following two solutions of the form (\ref{productcriticalpoint}) to (\ref{EulerLagrange})
$$\mathbb{S}^2(1)\times \mathbb{S}^1(\frac{1}{\sqrt{2}})\times \mathbb{S}^1(\frac{1}{\sqrt{2}}), \quad
\mathbb{S}^2(1)\times \mathbb{S}^1(\frac{1}{\sqrt{2}})\times \mathbb{S}^1(\frac{3}{\sqrt{10}}).$$
The value
$$\mathcal{L}_4[\mathbb{S}^2(1)\times \mathbb{S}^1(r_2)\times \mathbb{S}^1(r_3)]
=8\pi^3r_2r_3[\frac{7}{16}(4+x+y)^2-(8+x^2+y^2)],$$
In particular it is unbounded from above and below, and
$$\mathcal{L}_4[\mathbb{S}^2(1)\times \mathbb{S}^1(\frac{1}{\sqrt{2}})\times \mathbb{S}^1(\frac{1}{\sqrt{2}})]=48\pi^3, \quad
\mathcal{L}_4[\mathbb{S}^2(1)\times \mathbb{S}^1(\frac{1}{\sqrt{2}})\times \mathbb{S}^1(\frac{3}{\sqrt{10}})]=\frac{64\sqrt{5}}{3}\pi^3.$$

For $m=4$, that is $\mathbb{S}^1(r_1)\times \mathbb{S}^1(r_2)\times \mathbb{S}^1(r_3)\times \mathbb{S}^1(r_4)$,
let $x_i=\frac{1}{r_i^2}$, in the case we have
\begin{eqnarray*}
\mathcal{E}&=&\sum_{i=1}^4[2x_i^2-\frac{1}{2}(x_1^2+x_2^2+x_3^2+x_4^2)
+\frac{7}{8}(x_1+x_2+x_3+x_4)x_i
\\&&+\frac{7}{32}(x_1+x_2+x_3+x_4)^2]\frac{\nu_i}{r_i}.
\end{eqnarray*}
We have in the case only, modulo conformal equivalence, the following two solutions of the form (\ref{productcriticalpoint}) to (\ref{EulerLagrange})
$$\mathcal{L}_4[\mathbb{S}^1(1)\times \mathbb{S}^1(1)\times \mathbb{S}^1(1)\times \mathbb{S}^1(1)]=24\pi^4, \quad
\mathcal{L}_4[\mathbb{S}^1(1)\times \mathbb{S}^1(1)\times \mathbb{S}^1(1)\times \mathbb{S}^1(\frac{3}{\sqrt{5}})]=\frac{32\sqrt{5}}{3}\pi^4.$$
The value
$$\mathcal{L}_4[\mathbb{S}^1(r_1)\times \mathbb{S}^1(r_2)\times \mathbb{S}^1(r_3)\times \mathbb{S}^1(r_4)]=
8\pi^4r_1r_2r_3r_4[\frac{7}{16}(x_1+x_2+x_3+x_4)^2-(x_1^2+x_2^2+x_3^2+x_4^2)].$$
In particular it is unbounded from above and below, and
$$\mathcal{L}_4[\mathbb{S}^1(1)\times \mathbb{S}^1(1)\times \mathbb{S}^1(1)\times \mathbb{S}^1(1)]=24\pi^4, \quad
\mathcal{L}_4[\mathbb{S}^1(1)\times \mathbb{S}^1(1)\times \mathbb{S}^1(1)\times \mathbb{S}^1(\frac{3}{\sqrt{5}})]=\frac{32\sqrt{5}}{3}\pi^4.$$


\begin{thebibliography}{99}

\bibitem{Albin}
Albin, Pierre; Renormalizing curvature integrals on Poincare-Einstein manifolds. Adv. Math. 221 (2009), no. 1, 140--169.

\bibitem{AlexakisMazzeo}
Alexakis, Spyridon; Mazzeo, Rafe; Renormalized area and properly embedded minimal surfaces in hyperbolic 3-manifolds.
Comm. Math. Phys. 297 (2010), no. 3, 621--651.

\bibitem{Anderson82}
Anderson, Michael T. Complete minimal varieties in hyperbolic space. Invent. Math. 69 (1982), no. 3, 477--494.


\bibitem{Anderson83}
Anderson, Michael T. Complete minimal hypersurfaces in hyperbolic n-manifolds. Comment. Math. Helv. 58 (1983), no. 2, 264--290.

\bibitem{Anderson01}
Anderson, Michael T. $L^2$ curvature and volume renormalization of AHE metrics on 4-manifolds. Math. Res. Lett. 8 (2001), no. 1--2, 171--188.


\bibitem{ChangFangGraham}
Chang, Sun-Yung Alice; Fang, Hao; Graham, C. Robin; A note on renormalized volume functionals. Differential Geom. Appl. 33 (2014), suppl., 246--258.



\bibitem{ChangGurskyYang}
Chang, Sun-Yung A.; Gursky, Matthew J.;  An equation of Monge-Ampššre type in conformal geometry, and four-manifolds of positive Ricci curvature. Ann. of Math. (2) 155 (2002), no. 3, 709--787.

\bibitem{YangKingChang}
Chang, Sun-Yung A.; Qing, Jie; Yang, Paul C. A. Renormalized volumes for conformally compact Einstein manifolds. (Russian) ;
translated from Sovrem. Mat. Fundam. Napravl. 17 (2006), 129--142 J. Math. Sci. (N.Y.) 149 (2008), no. 6, 1755--1769

\bibitem{FeffermanGraham}
Fefferman, Charles; Graham, C. Robin; The ambient metric. Annals of Mathematics Studies, 178.
Princeton University Press, Princeton, NJ, 2012,  arXiv:0710.0919.

\bibitem{FeffermanGraham02}
Fefferman, Charles; Graham, C. Robin; Q-curvature and Poincare metrics. Math. Res. Lett. 9 (2002), no. 2--3, 139--151.

\bibitem{GlarosGoverHalbachWaldron}
M. Glaros, A. R. Gover, M. Halbasch and A. Waldron; Singular Yamabe Problem Willmore
Energies, arXiv:1508.01838.

\bibitem{Graham00}
Graham, C. Robin; Volume and area renormalizations for conformally compact Einstein metrics.
Rend. Circ. Mat. Palermo (2) Suppl. No. 63 (2000), 31--42.

\bibitem{Graham09}
Graham, C. Robin; Extended obstruction tensors and renormalized volume coefficients. Adv. Math. 220 (2009), no. 6, 1956--1985.

\bibitem{Graham17}
Graham, C. Robin; Volume renormalization for singular Yamabe metrics. Proc. Amer. Math. Soc. 145 (2017), no. 4, 1781--1792.

\bibitem{GrahamReichert}
Graham, C. Robin; Reichert, Nicholas; Higher-dimensional Willmore energies via minimal submanifold asymptotics.
arXiv:1704.03852

\bibitem{GrahamWitten}
Graham, C. Robin; Witten, Edward; Conformal anomaly of submanifold observables in AdS/CFT correspondence. Nuclear Phys. B 546 (1999), no. 1--2, 52--64.

\bibitem{GrahamZworski}
Graham, C. Robin; Zworski, Maciej; Scattering matrix in conformal geometry. Invent. Math. 152 (2003), no. 1, 89--118.

\bibitem{GoverWaldron1407}
A. R. Gover; A. Waldron; Submanifold conformal invariants and a boundary Yamabe
problem Conference on Geometrical Analysis-Extended Abstract, CRM Barcelona (2013),
arXived as Generalising the Willmore equation: submanifold conformal invariants from a
boundary Yamabe problem, arXiv:1407.6742.

\bibitem{GoverWaldron1506}
A. R. Gover; A. Waldron; Conformal hypersurface geometry via a boundary Loewner-
Nirenberg-Yamabe problem, arXiv:1506.02723.



\bibitem{GoverWaldron1603}
A. R. Gover; A. Waldron; Renormalized Volume, arXiv:1603.07367.


\bibitem{GoverWaldron1611}
A. R. Gover; A. Waldron; A calculus for conformal hypersurfaces and new higher
Willmore energy functionals, arXiv:1611.04055.

\bibitem{GoverWaldron1611B}
A. R. Gover; A. Waldron; Renormalized Volumes with Boundary, arXiv:1611.08345.


\bibitem{Gursky}
Gursky, Matthew J. The Weyl functional, de Rham cohomology, and K\"ahler-Einstein metrics. Ann. of Math. (2) 148 (1998), no. 1, 315--337.

\bibitem{Guven}
Guven, Jemal; Conformally invariant bending energy for hypersurfaces. J. Phys. A 38 (2005), no. 37, 7943--7955.


\bibitem{HanJiang}
Han, Qing; Jiang, Xumin; Boundary expansions for minimal graphs in the hyperbolic space.
arXiv: 1412.7608.

\bibitem{HanShenWang}
Han, Qing; Shen, Weiming; Wang, Yue; Optimal regularity of minimal graphs in the hyperbolic space.
Calc. Var. Partial Differential Equations 55 (2016), no. 1, Art. 3, 19 pp.

\bibitem{HardtLin}
Hardt, Robert; Lin, Fang-Hua; Regularity at infinity for area-minimizing hypersurfaces in hyperbolic space. Invent. Math. 88 (1987), no. 1, 217--224.



\bibitem{HenningsonSkenderis9806}
Henningson, M.; Skenderis, K. The holographic Weyl anomaly. J. High Energy Phys. 1998, no. 7, Paper 23, 12 pp.

\bibitem{HenningsonSkenderis9812}
Henningson, M.; Skenderis, K. Holography and the Weyl anomaly.
Quantum aspects of gauge theories, supersymmetry and unification (Corfu, 1998), 271--278,
Lecture Notes in Phys., 525, Springer, Berlin, 1999.

\bibitem{Lin89CPAM}
Lin, Fang-Hua; Asymptotic behavior of area-minimizing currents in hyperbolic space. Comm. Pure Appl. Math. 42 (1989), no. 3, 229--242.

\bibitem{Lin89Invent}
Lin, Fang-Hua; On the Dirichlet problem for minimal graphs in hyperbolic space. Invent. Math. 96 (1989), no. 3, 593--612.

\bibitem{Lin2012}
Lin, Fanghua; Erratum: On the Dirichlet problem for minimal graphs in hyperbolic space. Invent. Math. 187 (2012), no. 3, 755--757.




\bibitem{RyuTakayanagiPRL06}
Ryu, Shinsei; Takayanagi, Tadashi;
Holographic derivation of entanglement entropy from the anti-de Sitter space/conformal field theory correspondence.
Phys. Rev. Lett. 96 (2006), no. 18, 181602, 4 pp.

\bibitem{RyuTakayanagiJHEP06}
Ryu, Shinsei; Takayanagi, Tadashi; Aspects of holographic entanglement entropy. J. High Energy Phys. 2006, no. 8, 045, 48 pp.

\bibitem{Tonegawa}
Tonegawa, Yoshihiro; Existence and regularity of constant mean curvature hypersurfaces in hyperbolic space. Math. Z. 221 (1996), no. 4, 591--615.


\bibitem{Vyatkin}
Y. Vyatkin; Manufacturing conformal invariants of hypersurfaces, PhD thesis, University of Auckland, 2013.




\end{thebibliography}
\end{document}